\documentclass[9pto,reqno]{amsart}
\usepackage{cite}
\usepackage{array}
\usepackage{anysize}
\usepackage{fancyhdr}
\usepackage{float}
\usepackage{hyperref}
\usepackage{amssymb}
\usepackage{amsmath}
\usepackage{amsfonts}
\usepackage{amsthm}
\usepackage{mathrsfs} 
\usepackage{bbm}
\usepackage{bbold}
\usepackage{url}
\usepackage[section]{placeins}
\usepackage{graphicx}
\usepackage{subcaption}
\usepackage{caption}
\usepackage{color}
\usepackage{empheq}

\newtheorem{theorem}{Theorem}
\newtheorem{lemma}{Lemma}
\newtheorem{prop}{Proposition}
\newtheorem{Lemma}{Lemma}

\theoremstyle{definition}

\theoremstyle{remark}
\newtheorem{remark}{Remark}

\numberwithin{equation}{section}

\title[]{Pattern control via diffusion interaction}

\author[D. Ruiz-Balet]{Dom\`enec Ruiz-Balet}
\address[D. Ruiz-Balet]{Imperial College London, Department of Mathematics, Exhibition Rd, South Kensington, London SW7 2BX, United Kingdom}
\curraddr{}

\email{d.ruiz-i-balet@imperial.ac.uk}

\thanks{}

\author[E. Zuazua]{Enrique Zuazua}
\address[E. Zuazua]{Friedrich-Alexander-Universit\"at Erlangen-N\"urnberg, Department of Mathematics, Chair for Dynamics, Control, Machine Learning and Numerics (Alexander von Humboldt Professorship), Cauerstr. 11, 91058 Erlangen, Germany. 
	\newline \indent 
	Chair of Computational Mathematics, Fundaci\'on Deusto,	Avenida de las Universidades, 24, 48007 Bilbao, Basque Country, Spain. 
	\newline \indent
	Universidad Aut\'onoma de Madrid, Departamento de Matem\'aticas, Ciudad Universitaria de Cantoblanco, 28049 Madrid, Spain.}
\curraddr{}
\email{enrique.zuazua@fau.de}


\date{}

\dedicatory{}

\begin{document}
\maketitle

\begin{abstract}
We analyse a dynamic control problem for scalar reaction-diffusion equations, focusing on the emulation of pattern formation through the selection of appropriate active controls. While boundary controls alone prove inadequate for replicating the complex patterns seen in biological systems, particularly under natural point-wise constraints of the system state, their combination with the regulation of the diffusion coefficient enables the successful generation of such patterns.

Our study demonstrates that the set of steady-states is path-connected, facilitating the use of the staircase method. This approach allows any admissible initial configuration to evolve into any stationary pattern over a sufficiently long time while maintaining the system's natural bilateral constraints.

We provide also examples of complex patterns that steady-state configurations can adopt.

 \end{abstract}



\section{Introduction}

\subsection{Problem formulation and main results}\label{mot}
Reaction-diffusion equations are ubiquitous in modeling numerous natural phenomena, including patterns in embryos, species invasion, chemical reactions, and magnetic systems. Controlling these processes is crucial for various applications. In this article, we analyze the control of reaction-diffusion models to replicate some of the pattern formation phenomena observed in these complex systems.

We limit our analysis to the one-dimensional case, which, as we shall demonstrate, exhibits a rich and complex behavior that warrants careful, independent study.

Therefore, we consider the scalar time-evolving reaction-diffusion equation of the form:        \begin{equation}\label{diffusionparabolic}
     \begin{cases}
      \partial_t m-\partial_x\left(\mu(x,t)\partial_x m\right)=f(m)&\quad (x,t)\in (0,1)\times(0,T),\\
        m(0,t)=a_-(t); \, m(1,t)= a_+(t) &\quad t\in (0,T),\\
      m(x, 0)=m_0(x) &\quad x\in \{0,1\},
     \end{cases}
    \end{equation}
where $\mu$ is a positive diffusion control and $a$ is a boundary control with prescribed point-wise bounds 
\begin{equation}\label{controlconstraint}
-1\le a \le 1,
\end{equation}
 compatible with those imposed  on the state, 
 \begin{equation}\label{stateconstraint}
 -1\leq m\leq 1.
 \end{equation}

Bilateral state-constraints are natural in this context, since the solution describes the evolution of a density or a volumen fraction. Here, without loss of generality and to simplify the presentation we assume that $-1\leq m\leq 1$.

We will assume  that the diffusivity $\mu \in L^\infty((0,1)\times(0,T);\mathbb{R}^+)$ is a bounded measurable non-negative function depending both on space and time, the boundary control constituted by functions $$a_-\in L^\infty((0,T);[-1,1]); \,  a_+\in L^\infty((0,T);[-1,1])$$ corresponding to the controls at $x=0$ and $x=1$ respectively (satisfying the constraint \eqref{controlconstraint}), and the initial datum $m_0 \in L^\infty((0,1);[-1,1])$ (compatible with \eqref{controlconstraint}). This ensures the well-posedness of the system and the fulfilment of the state constraint \eqref{stateconstraint}.

One of the main distinguishing features of system \eqref{diffusionparabolic}, compared to those often considered in the literature, is the presence of two controls. The boundary control 
$a=a_\pm(t)$
 models an active control exerted by an external agent on the system's boundary. In contrast, the diffusion control 
$\mu(x,t)$ allows intervention in the medium's conditions and properties.

In our earlier papers, we analyzed whether the boundary control alone could steer the system from one equilibrium to another while preserving state constraints (\cite{loheac2017minimal,PoucholTrelatZuazua,Domenec-Zuazua,drift,loheac2021nonnegative}). We observed that, even for the simplest linear heat equation, a minimal control or waiting time emerges when aiming to reach the desired target under constrained dynamics. In the nonlinear setting, added barrier effects might appear, making some equilibria non-reachable. This largely depends on the domain's size, with the impact of boundary control weakening as the domain grows. Additionally, we proved that boundary control alone does not sufficiently modify the nature of these steady-states since the governing equations remain unchanged. This suggests the necessity of employing a second control, such as diffusivity, to regulate the dynamics effectively, as often seen in applications.

Our main objective in this article is to analyse the effect of combined strategies using both boundary and diffusion controls simultaneously. As we shall see, this combination allows us to shape the target equilibrium configurations and control the dynamics while preserving constraints.

Controlling diffusivity is natural in various applications. For instance,  in fluid mechanics, this can be achieved through temperature regulation \cite{einstein1906theory}. In ecology, diffusivity plays an essential role in determining species' survival and minimal area requirements \cite{PANDAS}. Moreover, in material science, the surface of metals can be patterned through diffusion control \cite{jiang2010diffusivity}. See subsection \ref{mot} for an in-depth discussion of modeling issues.

\begin{figure}
\captionsetup[subfigure]{justification=centering}
\centering
    \begin{subfigure}[c]{0.3\textwidth}
        \centering
         \includegraphics[width=\textwidth]{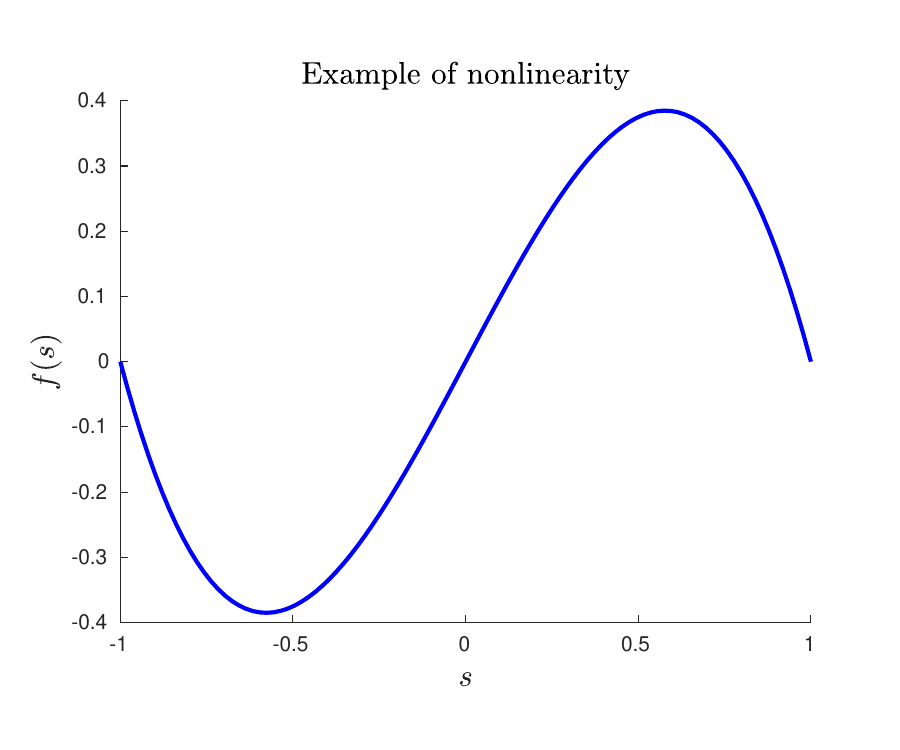}
        \caption{Bistable nonlinearity.}\label{F.NONLIN}
    \end{subfigure}
    \begin{subfigure}[c]{0.3\textwidth}
        \centering
         \includegraphics[width=\textwidth]{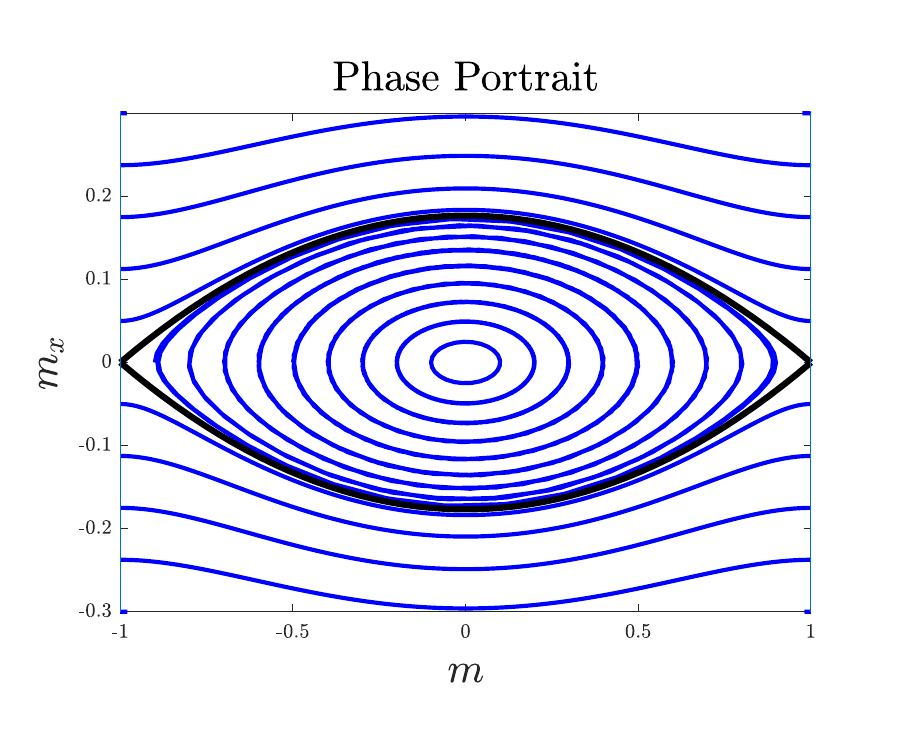}
        \caption{Phase portrait.}\label{F.PHASEPORTRAIT}
    \end{subfigure}
    \caption{Bistable non-linearity and typical phase portray of the ODE characterising equilibria.}\label{examples}
\end{figure}

Here and in the sequel $f\in C^1$ is a bistable nonlinearity satisfying  
\begin{equation}
f(-1)=f(0)=f(1)=0,  \, f'(-1),f'(1)<0, \,  f'(0)>0
\end{equation}
 as Figure \ref{F.NONLIN} shows.
We shall also assume that
\begin{equation}\label{notravel}
F(1)=\int_{-1}^1 f(s)ds=0
\end{equation}
which leads to the phase portraits as in  Figure \ref{F.PHASEPORTRAIT}. In the dynamic context the last condition guarantees, in particular, that travelling wave solutions behave like standing waves with null velocity of propagation, constituting steady-states.

We focus on the challenging problem of controllability, where the trajectory must be driven to a specified final target within a given time horizon. The infinite velocity of propagation in heat-like equations allows for the control of trajectories even within very short time horizons. However, standard procedures \cite{FC-ZZ,fczz2000} do not consider state constraints, which are essential in our applications. The oscillatory nature of the controls and controlled trajectories makes it impossible to fulfil the physical bilateral constraints unless the control time horizon is sufficiently long.

To understand what can be achieved for this parabolic model under state constraints, it is important to first analyse the structure of the set of steady-states, which are the most natural targets for the dynamics. These steady-states are characterized by the following elliptic system:
        \begin{equation}\label{introdif0}
     \begin{cases}
      -\partial_x\left(\mu(x)\partial_x m\right)=f(m)\qquad &x\in(0,1),\\
      m(0)=a_-; m(1)= a_+, \\
      -1 \le m \le + 1.
     \end{cases}
    \end{equation}
        
 We introduce the set
$$\mathcal{S}:=\{m(x) \text{ solution of \eqref{introdif0} with } \mu\in L^\infty((0,1);\mathbb{R}^+), a_\pm \in \mathbb{R} \},$$
constituted by the steady-states, solutions of  the system above when the diffusivity coefficient $\mu$ varies, together with the boundary values $a_\pm$.

We will show that this set of steady-states is path connected, a fact that plays an important role for the control of the evolution parabolic model since it  allows us to implement the  staircase method \cite{DARIO} (inspired by \cite{CORON-TRELAT}), that permits to track the evolution of the system along a path of steady-states.
This method has played a crucial role to develop  the theory of control of reaction-diffusion with point-wise constraints, see \cite{PoucholTrelatZuazua, Domenec-Zuazua, drift, lissy2021state}.

 In Theorem \ref{TH3} we prove that $S$
can approximate (with respect to the $L^2$-topology) a subset of the set of simple functions with a momentum restriction (see \eqref{consp}). In Theorem \ref{TH1} we present the main consequences for the control problem under consideration.

A straightforward change of variables in the steady-state system above allows us to transform variable diffusivity into a multiplicative potential on the nonlinearity. Similarly, in the parabolic problem, this transformation holds. The control problems under consideration also make sense in this case, where the potential multiplying the nonlinearity becomes the second control, complementing the effect of the boundary control. We will also consider this second model and control problem, as it appears in some relevant applications.

\subsection{Applications}\label{mot}
\textcolor{white}{.} In this subsection, we describe in more detail a number of models and relevant applications where the control of diffusivity (or the nonlinearity of the system) is both relevant and feasible.

Recent laboratory experiments \cite{JSHARPE} have demonstrated that limb formation can be intimately related to Turing instabilities \cite{turing1990chemical}, a counterintuitive process where a stable system, when combined with another stable diffusive process, can lead to unstable eigenmodes at high frequencies. In \cite{JSHARPE}, the number of fingers in embryos could be controlled by tuning the dynamics through spatial heterogeneities, resulting in unstable eigenfunctions that are spatially concentrated in the fingering areas. These experiments were the primary motivation for our work, where we aimed to construct and analyze the simplest model exhibiting related phenomena. 

In the context of interacting particle systems \cite{DEMASI} (see also \cite{durrett1995ten}), a rigorous micro-macro derivation is developed from an $N$-dimensional lattice $h\mathbb{Z}^N$, for $h>0$,  to get the  reaction-diffusion model in the limit as $h \to 0$. In this setting each $x\in h\mathbb{Z}^N$ represents a particle with spin $\sigma(x)=\pm 1$ subject to two types of interactions. The first one corresponds to the spin-exchange among neighbours, leading to the Laplacian when $h\to0$. The second one is a Glauber dynamics, a spontaneous spin-flip, depending on the spins of neighbouring particles, modelled by the nonlinearity. 

For ecological models of population dynamics \cite{barton2011spatial,FIFEBOOK} or spatial evolutionary games \cite{HOFBAUER,HOFBAUER2003evolutionary}, the reproduction and death processes are represented by nonlinearity, while diffusion accounts for spatial random motion. The size of the domain where the diffusion process occurs plays a crucial role in species survival \cite{PANDAS}, naturally leading to the question of constrained diffusivity control. Indeed, by scaling, the diffusivity coefficient and the size of the domain play reverse roles: the impact of boundary control decreases as the domain size increases or the diffusivity diminishes. In the nonlinear setting, small diffusivity or large domains enhance the existence of nontrivial steady-states. The combined effect of the two control actions plays an important role in this context. Multiplicative control  on nonlinearity corresponds to speeding up or slowing down the Glauber effects - reproduction and death of individuals in ecology- or the rate at which a game is played. Control in diffusion corresponds to limiting or enforcing species' movement between areas or slowing down or enhancing spin-exchange dynamics. These processes are well understood for the micro model \cite{DEMASI}.

Moreover, heterogeneities are intrinsic in nature, and their impact on system dynamics has been intensively studied \cite{cantrell1991effects,mazari2021fragmentation,mazari2020optimal}.

We also refer to  \cite{PoucholTrelatZuazua,Domenec-Zuazua,drift} for other related results in  control employing phase-plane analysis,  \cite{bressan2024controlled,bressan2022optimal} in the context of shaping traveling waves in relation to eradication of pests and  \cite{kobeissi2023tragedy} for explicit solutions to certain classes of mean-field games exhibiting the so-called tragedy of the commons.

\subsection{Structure}
\textcolor{white}{.}
The structure of the paper is the following:
\begin{enumerate}
\item In Section \ref{Core} we present the main results of the paper together with some of the fundamental aspects of the methodology we develop.
\item In Section \ref{SProofs} we prove the main results.

 \item In Section \ref{Open} we present some conclusions and  open problems, including the possible extensions to the multidimensional case.
\end{enumerate}

\section{Core phenomena and main results}\label{Core}
\subsection{The homogeneous case $\mu(x,t)=\mu$ and admissible paths}
We analyse the existence, nature  and role of admissible paths of steady-states in  the spatially  homogeneous and time-independent case where $\mu(x,t)=\mu$ is a constant, \cite{PoucholTrelatZuazua,Domenec-Zuazua,drift}.  Here  the viscosity $\mu>0$ plays the role of a design parameter, and not that of an active control. The viscosity being fixed, the control enters on the system solely through its boundary.

 Steady-state  solutions coincide with the restriction to $x \in (0, 1)$ of the solutions of the Hamiltonian dynamics:
\begin{equation}\label{ODE}
 \frac{d}{dx}\begin{pmatrix}
              m\\
              m_x
             \end{pmatrix}=
             \begin{pmatrix}
              m_x\\
              -\frac{1}{\mu}f(m)
             \end{pmatrix},
\end{equation}
which preserves the energy
\begin{equation}\label{energy} E(m,m_x)=\frac{1}{2}m_x^2+\frac{1}{\mu}F(m),
\end{equation}
where $$F(m)=\int_{-1}^mf(s)ds.$$

On the other hand, as mentioned above, in view of assumption \eqref{notravel}, the traveling wave solutions of
$$\partial_t m-\mu\partial_{xx}m=f(m),$$
are stationary. Therefore, travelling wave profiles are also solutions of \eqref{ODE}. 

The steady-states $m\equiv -1$ and $m\equiv 1$ correspond to the points $(-1,0)$ and $(1,0)$ in the phase plane, that are topological saddles for \eqref{ODE}. The state $m\equiv 0$ corresponds to $(0,0)$, which is a center (see figure \ref{F.PHASEPORTRAIT} for the representation of the phase-portrait). 
The traveling wave profiles are heteroclinic orbits  connecting $(-1,0)$ with $(1,0)$ and viceversa,  satisfying:
$$ 0=\frac{1}{2}m_x^2+\frac{1}{\mu}F(m),$$
denoted by the black curves in Figure \ref{F.PHASEPORTRAIT}.
Moreover, by the uniqueness of the solution of the Cauchy problem \eqref{ODE}, these curves define an invariant region in the phase-plane. Consequently, the solution of \eqref{ODE} for  any initial data   inside the invariant region  remains inside of it for all $x\geq 0$.

This invariance, together with the continuity of the solution of \eqref{ODE} with respect to the initial data, allows to build the paths of steady-states fulfilling the relevant constraints \eqref{stateconstraint}.
In fact, the set of admissible steady-states is path connected, as   the arguments below show.

Given two steady-states, $m_1$ and $m_2$,  inside the invariant region, we need to construct a curve inside the invariant region connecting the points
\begin{equation*}
  \begin{pmatrix}
   m_1(0)\\
   \frac{d}{dx}m_1(0)
  \end{pmatrix}
  \quad \text{ with }\quad
  \begin{pmatrix}
   m_2(0)\\
   \frac{d}{dx}m_2(0),
  \end{pmatrix}
\end{equation*}
 see \cite{PoucholTrelatZuazua,Domenec-Zuazua,drift}. The points of this curve constitute the initial data for \eqref{ODE}, leading to a path of admissible steady-states  linking $m_1$ to $m_2$. 

The constant diffusion $\mu$ being given, the staircase method allows to control the dynamics along these paths by means of the sole action of the boundary control. 

In particular, given that $(0,0)$ is always inside the invariant region, corresponding to the trivial steady-state $m\equiv 0$, the dynamics can be controlled from $m\equiv 0$ to any steady-state inside the invariant region, and vice-versa.

\subsection{Homogeneous time-dependent diffusion regulation, $\mu(x,t)=\mu(t)$}
\begin{figure}
\includegraphics[scale=0.35]{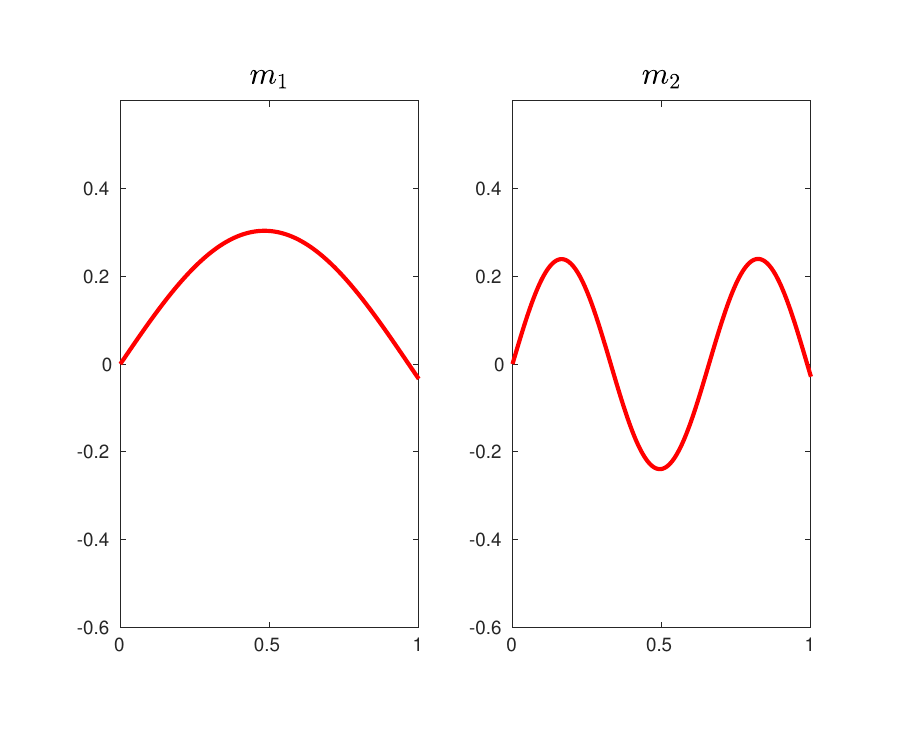}
 \caption{Steady-state solutions $m_1$ and $m_2$ with diffusivities $\mu_1>\mu_2$, constituted by arcs of a rescaled periodic orbit in the phase-plane.}\label{F.EXAMPLE}
\end{figure}

We now consider time-dependent  diffusion coefficients  $\mu$, i.e. $\mu(x,t)=\mu(t)$, playing the role of a second active control. In this case the equation \eqref{diffusionparabolic} reads
\begin{equation}\label{homogeneousregulation}
\begin{cases}
  \partial_t m-\mu(t)\partial_{xx}m=f(m)&\qquad (x,t)\in (0,1)\times (0,T),\\
 m=a(x,t)&\qquad (x,t)\in\{0,1\}\times (0,T),\\
 m(x,0)=m_0(x)\in L^\infty((0,1),[-1,1]).
\end{cases}
\end{equation}
By suitable choosing the time-dependent  diffusivity $\mu(t)$, which will now play the role of  a control function, one can control the system, for instance, from a positive function with  one single  maximum value $m_{\max}$ to a changing sing function withs two maxima of the same amplitude $m_{\max}$. See Figure \ref{F.EXAMPLE}, where the steady-states $m_i, i=1,2$ obey
\begin{equation*}
\begin{cases}
 -\mu_i\partial_{xx}m_i=f(m_i)&\qquad x\in(0,1),\\
 m_i=0&\qquad x\in\{0,1\},
  \end{cases}
\end{equation*}
for $i=1, 2$, with diffusivities $m_1>m_2$, so that $m_1$ has only one point of maximum, while $m_2$ has two of them. Note that, by diminishing the viscosity,  the frequency of oscillation of solutions increases, which allows for this construction.

The path of steady-states connecting $m_1$ with $m_2$, can be built patching two sub-arcs,  the first one connecting $m_1$ to the null sate, while the second one connects the null state  to $m_2$. More precisely, the following holds:

\begin{prop}\label{finger}
 Given two natural numbers, $k,n\in\mathbb{N}$ and  $0<\overline{m}<1$,  there exist two steady-states $m_1$ and $m_2$ for suitable constant diffusion coefficients $\mu_1,\mu_2\in\mathbb{R}^+$, with  $k$ and $n$  maxima within $x\in (0, 1)$, respectively,  taking the  value $\overline{m}$.
 
 These two patterns are path-connected in an admissible way. More precisely,   there exists $T^*_{\overline{m},n,k}>0$ such that for every $T\geq T^*_{\overline{m},n,k}$, there exist controls $\mu\in L^\infty((0,T),\mathbb{R}^+)$ and $a\in L^\infty((0,T),[-1,1]^2)$ such that the solution $m$ of \eqref{homogeneousregulation}, with initial datum $m_1$ at time $T$ satisfies that
 $$m(T;a,\mu,m_1)\equiv m_2.$$
 
\end{prop}

\begin{remark}
 The path of steady-states also solves
 $$ -\partial_{xx}m=\xi(t) f(m),$$
 with $\xi(t)=1/\mu(t)$, which can be interpreted as a multiplicative control entering in the nonlinearity of the system. 
 
 In other words, the class of steady-states is the same for both systems. Therefore, the  statement of Proposition \ref{finger} also applies for the multiplicative control
 $$ \partial_t m-\partial_{xx} m=\xi(t)f(m).$$
 The techniques of this paper can also be considered to consider the corresponding control problem:
   \begin{equation}\label{multiplicative}
     \begin{cases}
      \partial_t m-\partial_{xx}m=\xi(x,t)f(m)&\quad (x,t)\in (0,1)\times(0,T),\\
        m(0,t)=a_-(t); \, m(1,t)= a_+(t) &\quad t\in (0,T),\\
      m(x, 0)=m_0(x) &\quad x\in \{0,1\}.
     \end{cases}
    \end{equation}

\end{remark}

Once  the connectivity of the set of steady-states has been proved, the staircase method \cite{DARIO}, assures the controllability for large times from any initial data in $\mathcal{S}$  to any target data in $\mathcal{S}$ . 

\subsection{Heterogeneous diffusion $\mu(x,t)$}

Once the case of homogenous time-dependent diffusivity $\mu(t)$  has been addressed, we analyze that in which the diffusivity depends also in $x$, i.e. $\mu=\mu(x,t)$.

\subsubsection{Path connectivity}
\begin{figure}
 \includegraphics[scale=0.35]{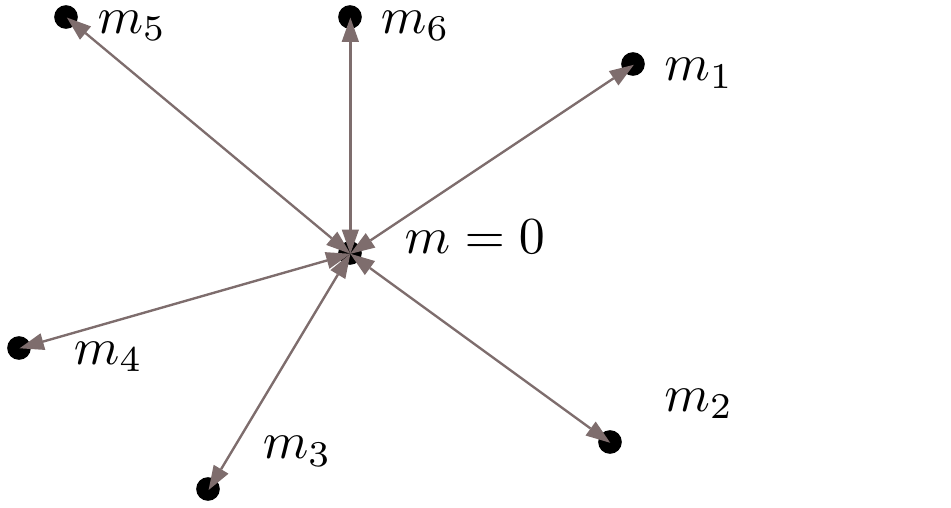}
 \caption{Star connectivity. All  admissible steady-states are connected by continuous paths of admissible steady-states with the trivial one $m=0$. In this figure,  six steady-states $m_1, ..., m_6$, represented with dots, are connected with $m=0$. 
 }\label{starconn}
\end{figure}

In this case too, all steady-states are path connected:
\begin{theorem}\label{TH0}
All elements of $\mathcal{S}$ are path connected in an admissible way. 
\end{theorem}
\begin{remark} As we shall see later on,  the connectivity of the set of steady-states and the staircase method assure the controllability for large times from any initial data in $\mathcal{S}$  to any target data in $\mathcal{S}$. 
\end{remark} 
To guarantee that all steady-states are path connected, it is sufficient to show that all of them can be connected to the trivial one, see Figure \ref{starconn}. Since the paths are reversible, the connectivity of any pair of steady-states will be automatically guaranteed. 

The connectivity towards the null steady-state will be proved by extending the domain of definition $x\in (0, 1)$, using that they correspond to arcs of trajectories of globally defined trajectories of the ODE.

Note also that for highly diffusive controls, the steady-state  $m=0$ is dynamically stable and attracts all admissible initial data for the parabolic problem. This, in turn, combined with the local controllability of the system by means of boundary control,  implies that any admissible initial datum  can be controlled to $m=0$ in finite time, and, from there, to any element of $\mathcal{S}$  following a path of steady-states. 

In this way we conclude that, departing from any initial configuration,  any state in  $\mathcal{S}$  can be reached in a sufficiently large finite time by means of suitable controls along admissible trajectories. Thus, to gain a complete understanding on the final patterns that the dynamics can reach, it suffices to describe accurately the set of equilibria $\mathcal{S}$.

To illustrate the nature of $\mathcal{S}$, hereafter we build a subset $\mathcal{S}^*$ of simple functions, with a momentum restriction, that can be approximated by elements of $\mathcal{S}$:
$$\mathcal{S}^*=\left\{ g\in L^\infty((0,1);(-1,1)): \quad g=\sum_{n=1}^{2N} \lambda_n\chi_{(a_{n-1},a_{n})}, \quad \text{with $a_0=0$, $a_{2N}=1$ and satisfying \eqref{consp}}\right\} $$

\begin{equation}\label{consp}
 \mathrm{sign}\left(\lambda_{n+1}\lambda_n\right)=-1,\quad |a_{n}-a_{n-1}||f(\lambda_n)|=|a_{n+1}-a_{n}| |f(\lambda_{n+1})|\quad n=1,3,...,2N-1.
\end{equation}

The equality $|a_{n}-a_{n-1}||f(\lambda_n)|=|a_{n+1}-a_{n}| |f(\lambda_{n+1})|$ is due to the transmission conditions at the diffusivity jumps and the changing sign of the parameters $\lambda$ is related to the nonlinear pendulum dynamics of the ODE under consideration. 

\begin{theorem}[Approximate steady-states for the diffusive interaction]\label{TH3}
For any function  $v\in \mathcal{S}^*$ and for all $\epsilon>0$ there exist $\mu_{\epsilon} \ge 0$ and $m_{\epsilon}\in L^\infty((0,1);[-1,1])$ solution of:
 \begin{equation*}
  \begin{cases}
   -\partial_{x}\left(\mu_\epsilon(x) \partial_x m_\epsilon\right)=f(m_\epsilon)&\qquad x\in(0,1)\\
   m_\epsilon=0&\qquad x\in\{0,1\}
  \end{cases}
 \end{equation*}
 such that 
  \begin{equation}\label{approximate}
 \|v-m_\epsilon\|_{L^2(0,1)}\leq \epsilon.
 \end{equation}
\end{theorem}

 Theorem \ref{TH3} will be proved constructing   piecewise constant diffusivities $\mu_\epsilon$, analysing the dynamics of the ODE characterising steady-states in the phase-plane. The resulting state $m_\epsilon$, even if continuous,  has  jump discontinuities in its derivative whenever $\mu$ has a discontinuity. Indeed, when $\mu$ has a discontinuity at $x_0$, the solution fulfils the following natural transmission conditions:
$$ \mu(x_0^-)=
            \mu_1, \,
            \mu(x_0^+)=\mu_2; \quad \lim_{x\to x_0^-} \mu_1\partial_x m(x)=\lim_{x\to x_0^+} \mu_2\partial_x m(x).$$

\begin{figure}
\centering
    \begin{subfigure}[c]{0.3\textwidth}
        \centering
 \includegraphics[width=\textwidth]{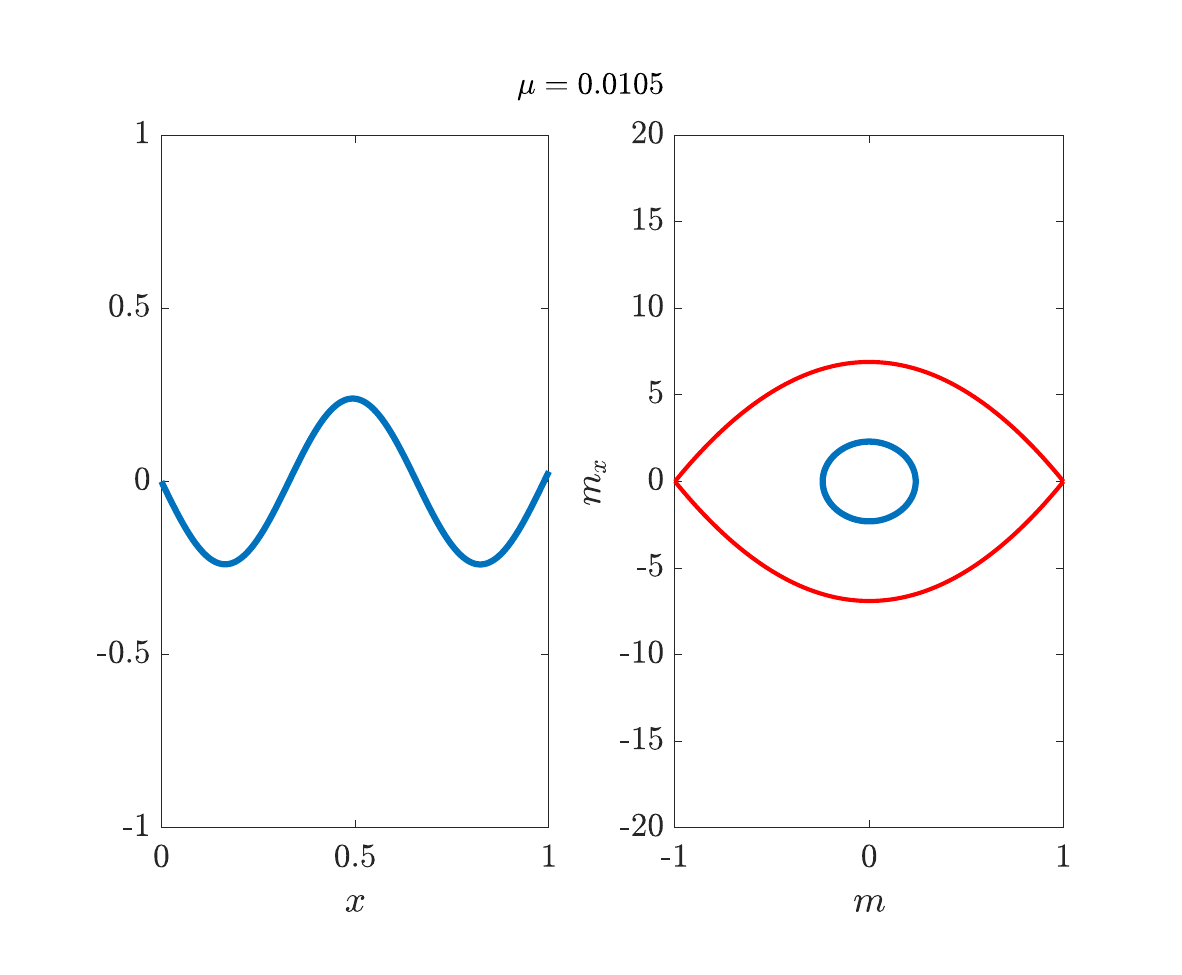}
        \caption{}
    \end{subfigure}
        \begin{subfigure}[c]{0.3\textwidth}
        \centering
 \includegraphics[width=\textwidth]{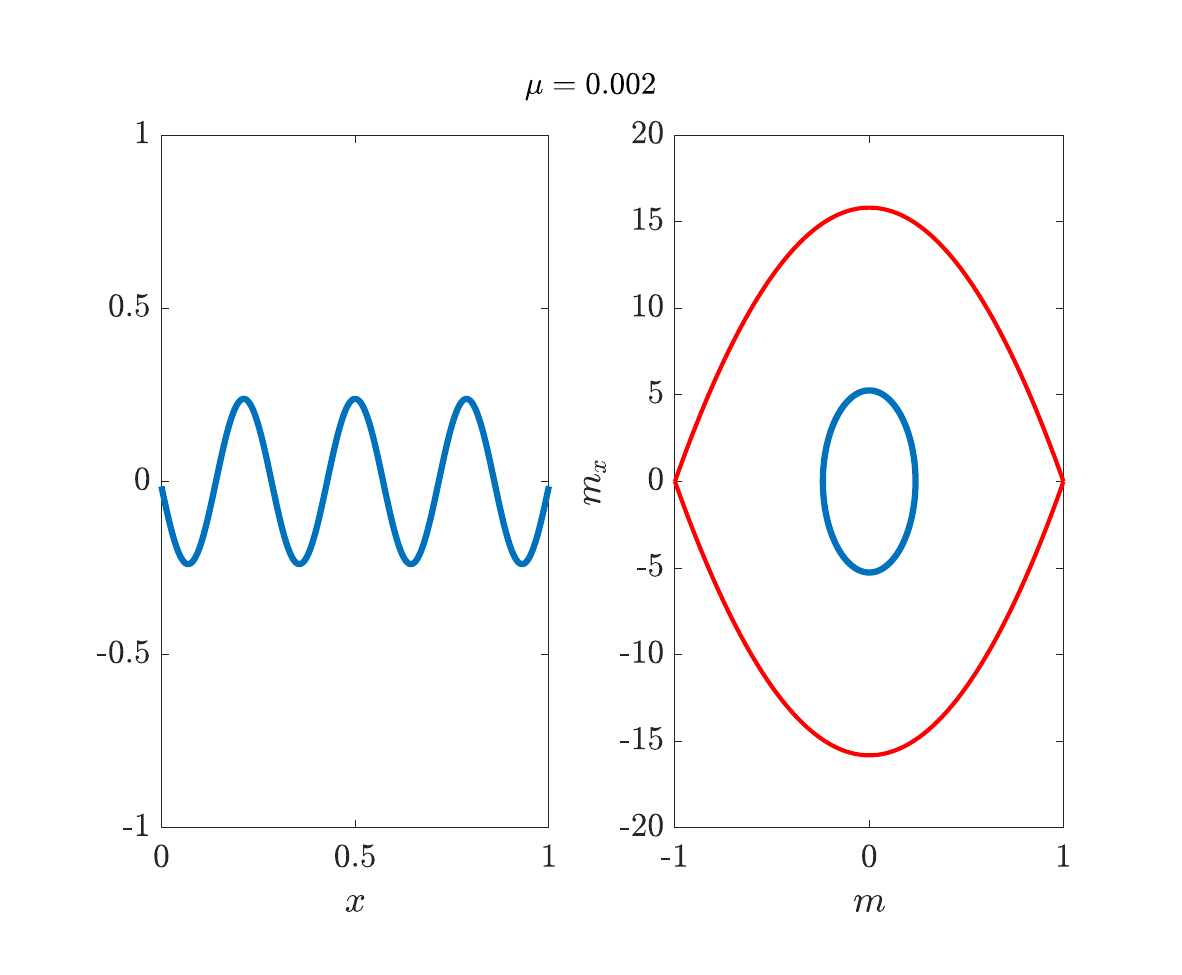}
        \caption{}
    \end{subfigure}
        \begin{subfigure}[c]{0.3\textwidth}
        \centering
 \includegraphics[width=\textwidth]{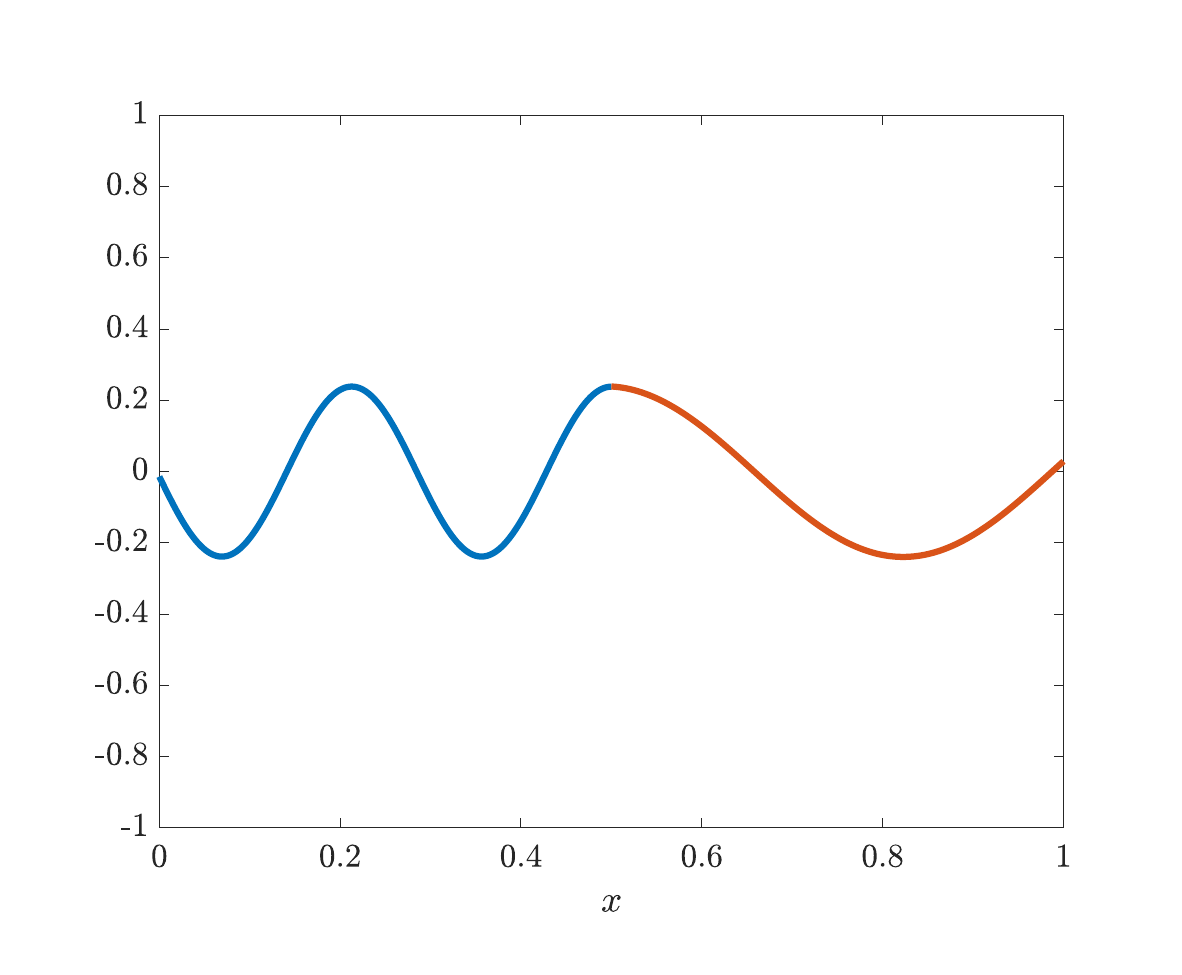}
        \caption{}
    \end{subfigure}
\caption{Example of a composite pattern  within $\mathcal{S}$. (A) Steady-state for a constant diffusivity $\mu_1$ with a critical point at $x=0.5$; (B)  Steady-state for a constant diffusivity $\mu_2<\mu_1$ with a critical point at $x=0.5$; (C) Steady-state presenting a frequency contrast, with an heterogeneous diffusivity combining the previous patterns in different regions.}\label{difconex}
\end{figure}

The set $\mathcal{S}^*$ does not exhaust the class of patterns that can be approximated through elements of $\mathcal{S}$, as shown  in Proposition \ref{finger}. In fact, as illustrated in Figure \ref{difconex}, one can also build continuous composite profiles matching two different oscillatory patterns.

\subsection{Some model variants}
The same questions arise and our techniques can be adapted for the following model variants
 \begin{equation*}
  \begin{cases}
   \partial_t m-\mu(x,t)\partial_{xx} m=f(m)&\qquad x\in(0,1),\\
   m(0)=a_-(t)\in [-1,1], \quad m(1)=a_+(t)\in [-1,1]&
  \end{cases}
 \end{equation*}
where the control enters multiplying the Laplacian or  \begin{equation}\label{eq.multiplicative}
  \begin{cases}
   \partial_t m-\partial_{xx} m=\xi(x,t)f(m)&\qquad x\in(0,1),\\
   m(0)=a_-(t)\in [-1,1], \quad m(1)=a_+(t)\in [-1,1]&
  \end{cases}
 \end{equation}
where it enters multiplying the nonlinearity.

This type of control falls into the class of the multiplicative controls studied by \cite{cannarsa2017multiplicative} and others.

The steady states of \eqref{eq.multiplicative} fulfil
 \begin{equation}\label{eq.multiplicative}
  \begin{cases}
   -\partial_{xx} m=\xi(x)f(m)&\qquad x\in(0,1),\\
   m(0)=a_-\in [-1,1], \quad m(1)=a_+\in [-1,1]&
  \end{cases}
 \end{equation}

For this model the set of controlled steady-states turns out to be dense:

\begin{theorem}\label{TH2}
 For any function  $v\in L^\infty((0,1);[-1,1])$ and for all $\epsilon>0$ there exists $\xi\in L^\infty((0,1);\mathbb{R}_+)$ such that $m_{\epsilon}\in L^\infty((0,1);(-1,1))$ solution of:
 \begin{equation*}
  \begin{cases}
   -\partial_{xx} m_\epsilon=\xi(x)f(m_\epsilon)&\qquad x\in(0,1),\\
   m_\epsilon=0&\qquad x\in\{0,1\},
  \end{cases}
 \end{equation*}
 such that $\|v-m_\epsilon\|_{L^2((0,1))}\leq \epsilon$.
\end{theorem}
The proof is very similar to the proof of Theorem \ref{TH3}. In consists in approximating the target state by a piecewise-constant function and building the control to approximate the latter.

For the present model the transmission condition is not required and this allows for the density of the set of steady-states.

Theorem \ref{TH0} also holds for the set of steady-states of \eqref{eq.multiplicative}.

Theorem \ref{TH2} can be proved analysing the affine ODE control problem
\begin{equation}\label{ODEcontr}
 \frac{d}{dx}\begin{pmatrix}
              m\\
              m_x
             \end{pmatrix}=\begin{pmatrix}
             m_x\\
             -\xi(x)f(m)
             \end{pmatrix},
\end{equation}
where $\xi\in L^\infty((0,1);\mathbb{R}^+)$ is the control function. 

\begin{figure}
\captionsetup[subfigure]{justification=centering}
\centering
    \begin{subfigure}[c]{0.225\textwidth}
        \centering
         \includegraphics[width=\textwidth]{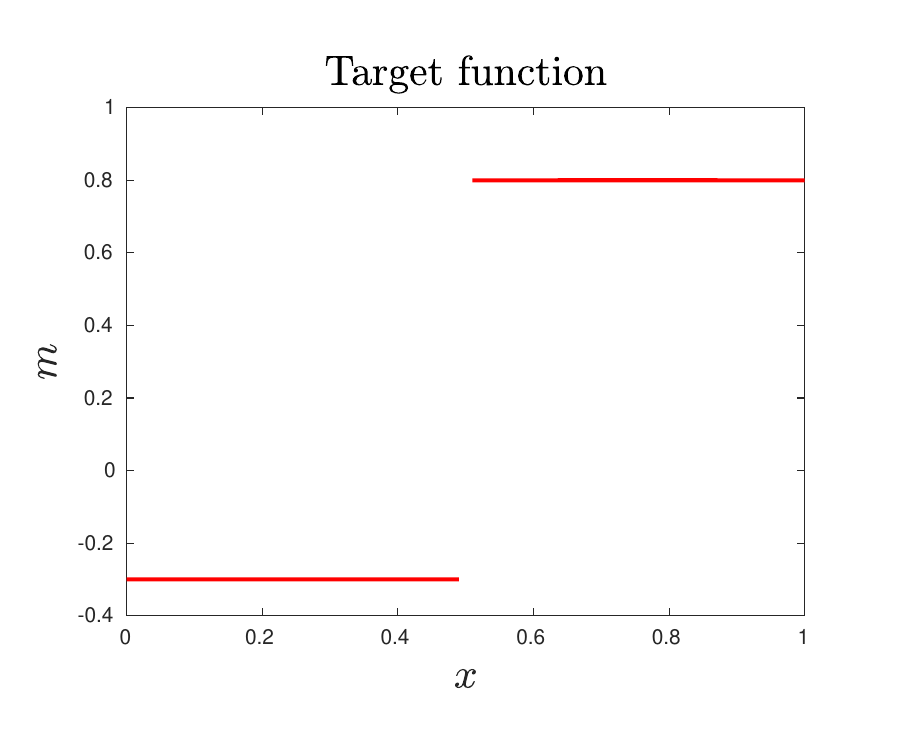}
        \caption{}
    \end{subfigure}
    \begin{subfigure}[c]{0.225\textwidth}
        \centering
         \includegraphics[width=\textwidth]{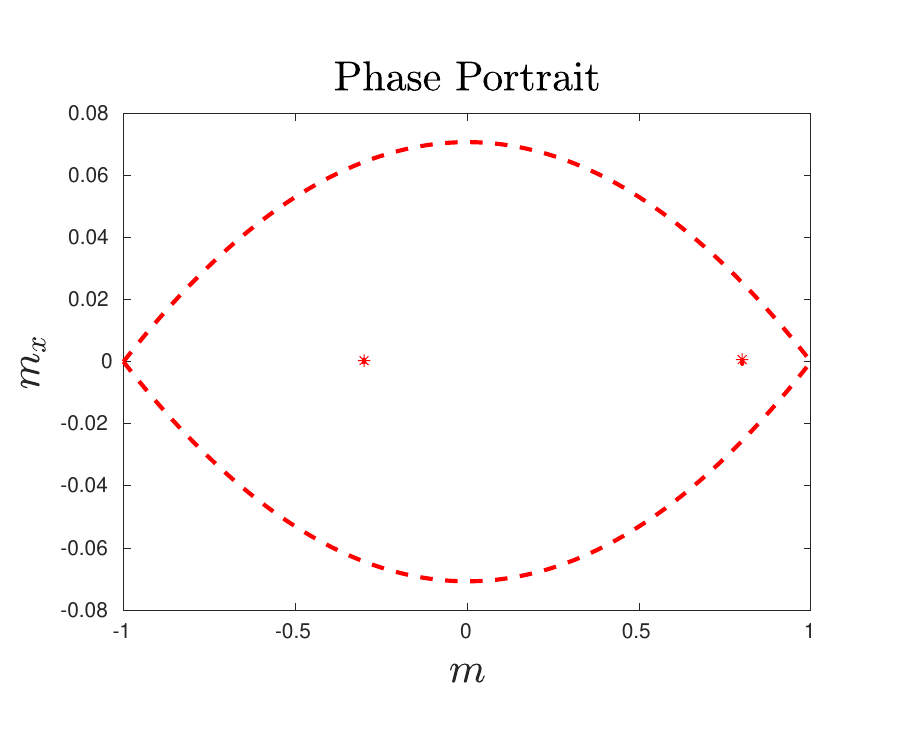}
        \caption{}
    \end{subfigure}
        \begin{subfigure}[c]{0.225\textwidth}
        \centering
 \includegraphics[width=\textwidth]{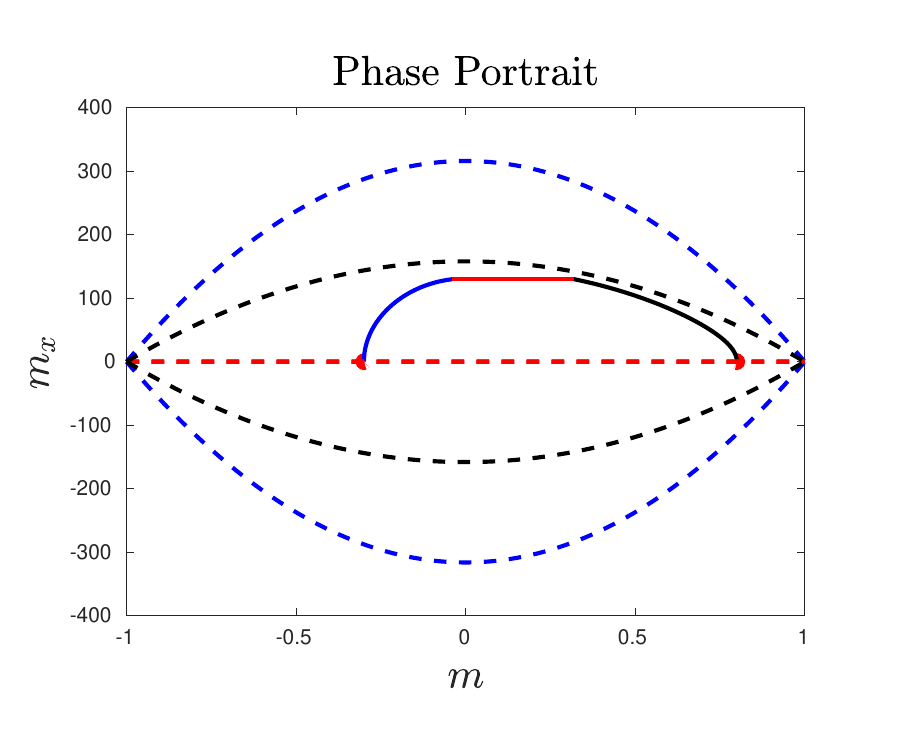}
        \caption{}
    \end{subfigure}
    \begin{subfigure}[c]{0.225\textwidth}
        \centering
 \includegraphics[width=\textwidth]{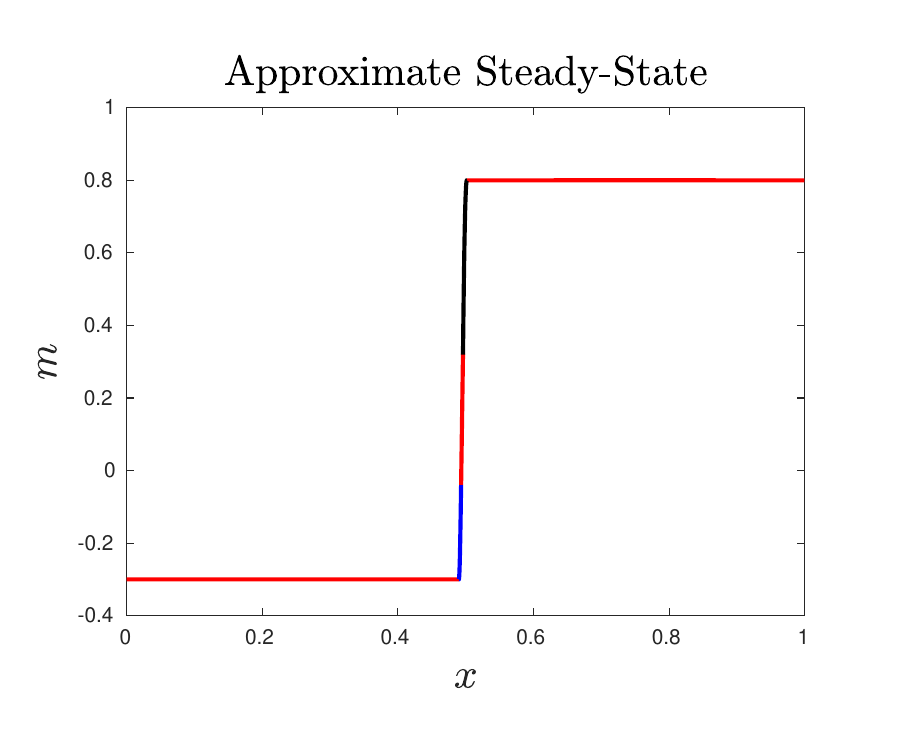}
        \caption{}
    \end{subfigure}
 \caption{Approximation using the affine controllability of the ODE system. (A) A target step function (B) its representation in the phase plane as two points, the dashed red line indicates the invariant region enclosed by the standing traveling waves for a certain constant multiplicative action  (C)  Controllability of the associated ODE system driving the state from one point to the other employing different values of the multiplicative control (multiple colours). The different dashed lines indicate the invariant region delimited by the corresponding standing traveling waves with the associated value of the multiplicative control. (D) The corresponding state in the physical space.}\label{controlaffineex}
\end{figure}

Recall that, by changing $\xi$,  the invariant region changes. This can  help generating the steady states we desire. See Figure \ref{controlaffineex} for an illustration of how the multiplicative control can be used for approximation.

The functions $\xi$ we employ are piecewise constant in $x$. By smoothly changing its switching points one obtains  the connected family of steady-states allowing to apply the staircase method.

A key point for proving Theorem \ref{TH2} is the following Lemma, which guarantees the exact controllability of \eqref{ODEcontr} provided the initial data and  target are in the region
$$R:=\left\{(m,m_x)\in\mathbb{R}^2:\quad -1<m<1\right\}.$$

\begin{lemma}\label{Col1}
 For any $L>0$,  any initial data $\boldsymbol{m_0}\in R\backslash\{(0,0)\}$, and any target $\boldsymbol{m_L}\in R\backslash\{(0,0)\}$, there exists  a strictly positive function $\xi\in L^\infty((0,L),\mathbb{R}^+)$ such that the solution of \eqref{ODEcontr} with initial data $\boldsymbol{m_0}$ takes the value $\boldsymbol{m}(L;\xi,\boldsymbol{m_0})=\boldsymbol{m_L}$ at $x=L$, and,  furthermore, $\boldsymbol{m}(x;\xi,\boldsymbol{m_0})\in R$ for all $0\leq x\leq L$.
\end{lemma}

\subsection{Dynamic control}

Using the results above on the nature  of the set $\mathcal{S}$  and the staircase method the following control result can be proved:

\begin{theorem}[Approximate controllability with point-wise constraints]\label{TH1}
 Take any pair of functions $m_0\in L^\infty ((0,1);[-1,1])$ and $m_1 \in \mathcal{S}$. Then, for every $\epsilon>0$ there exists $T_{\epsilon,m_1}>0$ such that for all $T\geq T_{\epsilon,m_1}>0$, there exists $\mu_\epsilon\in L^\infty((0,1)\times(0,T);\mathbb{R}^+)$, $a_\epsilon\in L^\infty(\{0,1\}\times(0,T);[-1,1])$ such that the solution  of \eqref{diffusionparabolic} satisfies $m(t;a_\epsilon,\mu_\epsilon,m_0)\in L^\infty((0,1);[-1,1])$ a.e. in $(0,T)$ 
and
\begin{equation*}
 \|m(T;a_\epsilon,\mu_\epsilon,m_0)-m_1\|_{L^2}\leq \epsilon.
\end{equation*}
Moreover $m(T; a_\epsilon,\mu_\epsilon,m_0)$ is a steady-state in $\mathcal{S}$, with diffusivity $\mu_\epsilon(x,T)$ and boundary values $a_\epsilon(x,T)$ at $x=0, 1$.
\end{theorem}

\begin{remark}
In fact, the control time $T_{\epsilon,m_1}$   can be taken to depend on the target $m_1$ and $\varepsilon$, but  independent of $m_0$. This is due to the exponential attractiveness of the null state of the system \eqref{diffusionparabolic},  that is enhanced  when the diffusivity $\mu$ increases. More precisely, the convergence towards the null equilibrium of the solutions of \eqref{diffusionparabolic} as $t\to \infty$ is faster than the one determined by the eigenvalue 
$$\lambda_1(\mu):=\min_{v\in H^1_0}\frac{\int_0^1\mu|\nabla v|^2-\|f'\|_\infty v^2dx}{\int_0^1 v^2dx},$$
which, obviously, fulfils
$$\lim_{\mu\to +\infty}\lambda_1(\mu)=+\infty.$$
 Therefore, by means of a suitable choice of the diffusivity one may stabilise the system towards $m\equiv 0$ arbitrarily fast and  then apply the staircase method to drive the system towards the final state $m_1$ (see \cite{DARIO}). This second regime in which the system is controlled towards $m_1$ up to an $\varepsilon$ error needs of a control time depending both on $m_1$ and $\varepsilon$. But the overall control time turns out to be independent of $m_0$.
\end{remark}

\section{Proofs}\label{SProofs}

\subsection{Proof of Proposition \ref{finger}}\label{prooffinger}
It suffices to set the diffusivity as
\begin{equation*}
 \mu(s)=\begin{cases}
         \mu_1 \quad s\in[0,1/2),\\
         \mu_2 \quad s\in[1/2,1],
        \end{cases}
\end{equation*}
and then define the boundary values $a(s)$ so that
\begin{equation*}
 a(s)=\begin{cases}
         \text{A path connecting $m_1$ to $m\equiv 0$ for } s\in[0,1/2),\\
         \text{A path connecting $m\equiv 0$ to $m_2$ for }  s\in[1/2,1].
        \end{cases}
\end{equation*}
This is one of the many possible paths connecting two elements of $\mathcal{S}$. \begin{figure}
\includegraphics[scale=0.35]{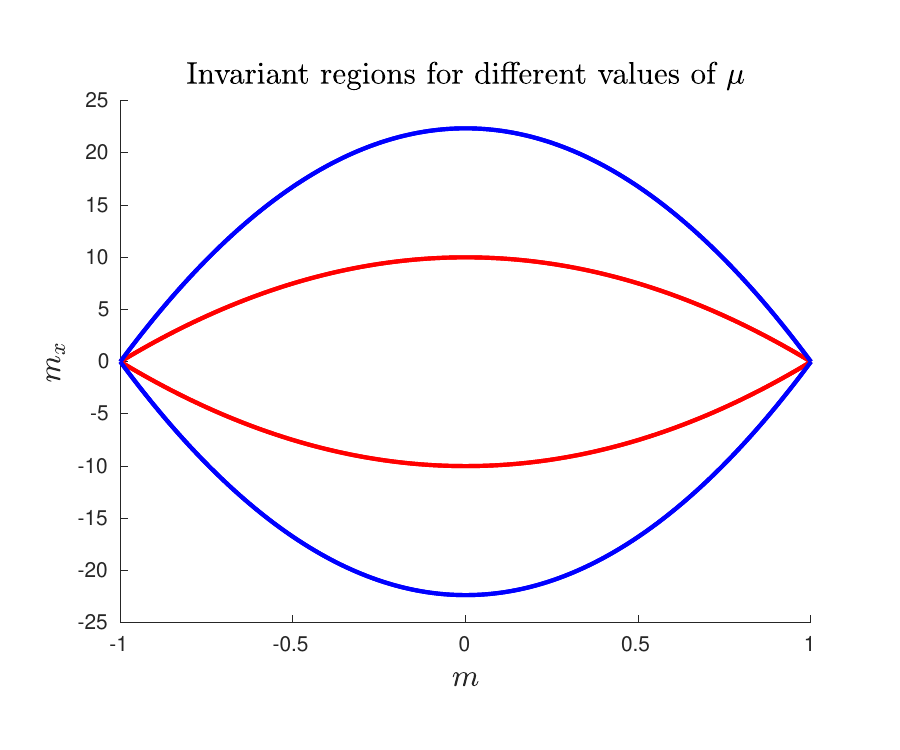}
 \caption{Invariant regions for different values of the diffusivity $\mu$ for the symmetric nonlinearity $f(s)=s(s+1)(s-1)$. In blue $\mu=0.001$ in red $\mu=0.005$.}\label{F.INVARIANTREGIONCHANGE}
\end{figure}

\subsection{Preliminary results}
The ideas of the proofs of the Theorems
  are similar: switching from one control to another and combining the use of periodic orbits of different frequencies. 
  
  Let us first analyse how the period behaves when the diffusion control is constant and $\mu\to 0$ (resp. $\mu\to +\infty$).

 Let $\Gamma_\mu$ be the invariant region defined as
 \begin{equation*}
 	\Gamma_\mu:=\left\{(m,m_x)\in [-1,1]\times \mathbb{R}:\quad m_x^2\leq \frac{-2}{\mu}F(m)\right\} 
 \end{equation*}

\begin{Lemma}\label{Period}
 \begin{enumerate}
 \item Fix $m_{\min}\in(-1,0)$ (or $m_{\max}\in(0,1)$). The period of the periodic orbit starting at $(m_{\min},0)$ (or $(m_{\max},0)$) tends to $0$ as $\mu\to 0$.
  Moreover, fixing the maximum value of an orbit, the period  of the orbit passing through that maximum depends continuously on  $\mu$.
  
  \item For every $a,b\in [-1,1]$ and any $L$ there exists $\mu\in\mathbb{R}_+$ and $(m,m_x)\notin \Gamma_\mu$ such that the solution of \eqref{ODE} starting at $(m,m_x)$ satisfy that $m(L^*)=b$ with $0<L^*\leq L$. Moreover, $L^*\to0$ when $\mu\to +\infty$.

  \end{enumerate}

\end{Lemma}
\begin{proof}
We proceed in order.
\begin{enumerate}
\item 
The system is symmetric with respect to the $m_x=0$ axis. Let $\tau$ be the period of one orbit inside the invariant region defined by the static traveling waves. The period $\tau$ can be computed as follows:
Consider an arc of an orbit that goes from  the minimum value $m_{\min}$ of the orbit to the maximum $m_{\max}$. The map $m:[0,\tau/2]\to [m_{\min},m_{\max}]$ is invertible, since both the positive and the negative branch of the curve $m_x(m)$  only change sign in $\{m_{\min},m_{\max}\}$. Hence $m_x$ has the same sign along $[0,\tau/2]$.
 \begin{align*}
 \frac{\tau}{2}&=\int_0^{\frac{\tau}{2}}dx=\int_{m^{-1}(m_{\min})}^{m^{-1}(m_{\max})}dx=\int_{m_{\min}}^{m_{\max}}\frac{1}{m_x(m)}dm\\
  &=\int_{m_{\min}}^{m_{\max}} \frac{dm}{\sqrt{2\left(E-\frac{1}{\mu}F(m)\right)}}=\sqrt{\mu}\int_{m_{\min}}^{m_{\max}} \frac{dm}{\sqrt{2\left(\mu E-F(m)\right)}}.
 \end{align*}
Note that the expression above is continuous with respect to $\mu$.
 
If $m_{\min}$ and $m_{\max}$ are away of $-1$ and $1$ then the function $G(m)=\mu E-F(m)$ near $m=m_{\min}$ behaves as $-f(m_{\min})(m-m_{\min})+\mathcal{O}(m^2)$. Similarly one can see that $G(m)$ behaves linearly around $m_{\max}$ which enables us to assure that
$$\int_{m_{\min}}^{m_{\max}} \frac{dm}{\sqrt{2\left(\mu E-F(m)\right)}}\leq C<+\infty,$$
where $C$ is independent of $\mu$, since the energy of the orbit $E$ is $E=\frac{1}{\mu}F(m_{\min})=\frac{1}{\mu}F(m_{\min})$. 

Note that any extreme of the orbit can take the value $0$ since, at that point, we would have $m_x=0$, moreover $m=0$ is a zero for $f$, and $(0,0)$ is a critical point for the ODE system.

The maximum value on the $m_x$ variable that trajectories in the invariant region reach is
\begin{equation}\label{maxmx}
\frac{1}{\mu}\sqrt{2(\mu E-F(0))}, 
\end{equation}
 which goes to infinity as $\mu$ goes to zero. Moreover, the derivative with respect to $m$ of the positive branch of $m_x(m)$ is
$$\frac{d}{dm}m_x=\frac{1}{\sqrt{\mu}}\frac{f(m)}{\sqrt{2(\mu E-F(m))}}$$
and goes to infinity as $\mu\to 0$ if $m\neq 0$.

This shows that small diffusivities produce a significant increase in the gradient in a very short interval. 
\item 
Let us now analyze what happens when the gradient is large and  with high diffusivity $\mu$, for example, for trajectories outside the invariant region. They take values from $m=-\infty$ to $m=+\infty$ (or vice versa). Thus, following  trajectories outside the invariant region, we can access every value of $m$. 

When $\mu>0$ is large enough the space $x$ required to go from the value $m=a$ to the value $m=b$ with $b>a$ and $m_x(0)>0$ is approximately $(b-a)/m_x(0)$
$$x=\int_a^b \frac{dm}{\sqrt{m_x(0)^2+\frac{2}{\mu}(F(m(0)-F(m)))}}\to\frac{b-a}{m_x(0)}$$
when $\mu\to \infty$.

Then, for $\mu$ large and $m_x(0)$ large enough we can access to any $b>a$  in a small space interval $[0,x]$.
\end{enumerate}
\end{proof}

\subsection{Proof of Theorem \ref{TH3}}

\begin{enumerate}
 \item Let $g$ be a piecewise constant function of the form
\begin{equation}\label{ggg}
 g=\sum_{i=1}^N \alpha_i \chi_{(a_{i-1},a_i)}
\end{equation}
with $a_0=0$, $a_N=1$ and $a_{i}<a_{i+1}$.

We now explore the class of $\alpha$'s that can be approximated by steady-state solutions.

First of all we approximate each constant segment by an elliptic equation. Consider $I_i=(a_{i-1}+\delta,a_i-\delta)$ and the elliptic equation:
\begin{equation}\label{approxPB}
 \begin{cases}
  -\mu_i\partial_{xx}m_i=f(m_i)\quad &x\in I_i\\
  m_i=\alpha_i\quad &x\in\partial I_i
 \end{cases}
\end{equation}
Since $\alpha_i\in (-1,0)\cup (0,1)$, equation \eqref{approxPB} has a solution taking values in $(-1,1)$. To obtain 
$$ \|m_i-\alpha_i\|_{L^\infty(I_i)}\leq\epsilon$$
we can fix the maximum or minimum of \eqref{approxPB} to be at $\alpha_i\pm \epsilon$. Since the energy \eqref{energy} of the ODE associated to the problem \eqref{approxPB}
is preserved, we have a formula for the Neumann trace of a solution of \eqref{approxPB}:
\begin{equation}\tag{A}\label{EnergyTrace}
 \partial_x m_i^2=\frac{2}{\mu_i}\left(F(\alpha_i\pm\epsilon)-F(\alpha_i)\right).
\end{equation}
\item 
The derivative of $dp/dx$, with $p=\mu(x)\partial_xm$, is bounded regardless of the control action. In particular, if we are trying to approximate a piecewise constant function, since $p$ is continuous and its derivative is bounded, we can only expect to approximate piecewise constant functions fulfilling the transmission condition. 

Assume that $\mu$ is piecewise constant with a discontinuity at $x_0\in(0,1)$. Then any solution of
\begin{equation*}
  \begin{cases}
  -\partial_{x}\left(\mu(x)\partial_xm\right)=f(m)\quad &x\in (0,1)\\
  m=0\quad &x\in\{0,1\}
 \end{cases}
\end{equation*}
satisfies 
\begin{equation}\tag{B}\label{Tr}
\lim_{x\to x_0^-}\mu(x) \partial_x m(x)=\mu_1 \partial_x m_1= \mu_2 \partial_x m_2=\lim_{x\to x_0^+}\mu(x) \partial_x m(x)\end{equation}
\item Let $L_i$ be the length of the interval $I_i$. As in Lemma \ref{Period}, we have an expression  relating the diffusivity $\mu_i$, the energy of the system and $\alpha_i$ with $L_i$:
\begin{equation}\tag{C}\label{Length}
 L_i^2= 2\mu_i \left(\int_{\alpha_i\pm \epsilon}^{\alpha_i} \frac{dm}{\sqrt{F(\alpha_i\pm \epsilon)-F(\alpha_i)}}\right)^2
\end{equation}
\item Combining \eqref{EnergyTrace},\eqref{Tr} and \eqref{Length} one arrives at the following relationship between the $\alpha$'s depending on $\epsilon$ and $L_i$'s for two adjacent constant pieces of \eqref{ggg}:
\begin{equation}\label{desastre}
\left(\frac{\int_{\alpha_i\pm \epsilon}^{\alpha_i} \frac{dm}{\sqrt{F(\alpha_i\pm \epsilon)-F(\alpha_i)}}}{\int_{\alpha_{i-1}\mp \epsilon}^{\alpha_{i-1}} \frac{dm}{\sqrt{F(\alpha_{i-1}\mp \epsilon)-F(\alpha_{i-1})}}}\right)^{-2}=\frac{L_{i-1}^2}{L_i^2}\frac{F(\alpha_{i-1}\mp\epsilon)-F(\alpha_{i-1})}{F(\alpha_i\pm\epsilon)-F(\alpha_i)}
\end{equation}
Now we will see that when $\epsilon \to 0$ it corresponds to require that $g$ in \eqref{ggg} satisfies
\begin{equation}\label{conclusionR}
 L_{i-1}^2f(\alpha_{i-1})^2=L_{i}^2f(\alpha_i)^2.
\end{equation}
The right hand side of \eqref{desastre}, when $\epsilon\to 0$ tends to
$$ \frac{L_{i-1}^2}{L_i^2}\frac{f(\alpha_{i-1})}{f(\alpha_i)}$$
While for the left hand side we need to analize the integral
\begin{align*}
\int_\alpha^{\alpha+\epsilon}\frac{dm}{\sqrt{F(\alpha+\epsilon)-F(m)}}&=\int_0^\epsilon \frac{ds}{\sqrt{F(\alpha+\epsilon)-F(\alpha+s)}}\\
&=\int_0^\epsilon \frac{ds}{\sqrt{(\epsilon-s)(f(\alpha))+\mathcal{O}(\epsilon^2.s^2)}}\\
&\sim \frac{2\sqrt{\epsilon}}{\sqrt{f(\alpha)}}.
\end{align*}
Therefore, the left hand side of \eqref{desastre} tends to
\begin{equation*}
 \frac{f(\alpha_i)}{f(\alpha_{i-1})}.
\end{equation*}
Thus, one arrives to \eqref{conclusionR}.
\end{enumerate}

Consider a function of the form:
$$S=\left\{ g=\sum_{n=1}^{2N} \lambda_n\chi_{I_n},\quad \mathrm{sign}\left(\lambda_{n+1}\lambda_n\right)=-1,\quad L_n|f(\lambda_n)|=L_{n+1}|f(\lambda_{n+1})|\quad n=1,3,...,2N-1\right\} $$
\begin{equation}
 \sum_{n=1}^{2N} \lambda_n\chi_{I_n}
\end{equation}
with 
\begin{equation}
 \bigcup_{n=1}^{2N} I_n=(0,1)\qquad I_n\cap I_{n'}=\emptyset\quad\text{ if }n\neq n'.
\end{equation}
Let $L_n=|I_n|$.

If 
\begin{equation}
 \mathrm{sign}\left(\lambda_{n+1}\lambda_n\right)=-1,\quad L_n|f(\lambda_n)|=L_{n+1}|f(\lambda_{n+1})|\quad n=1,3,...,2N-1
\end{equation}
 there exists $\mu\in L^\infty([0,1];\mathbb{R}^+)$ such that \eqref{approximate} holds with $v=g$.

\begin{figure}
\captionsetup[subfigure]{justification=centering}
\centering
    \begin{subfigure}[c]{0.3\textwidth}
        \centering
         \includegraphics[width=\textwidth]{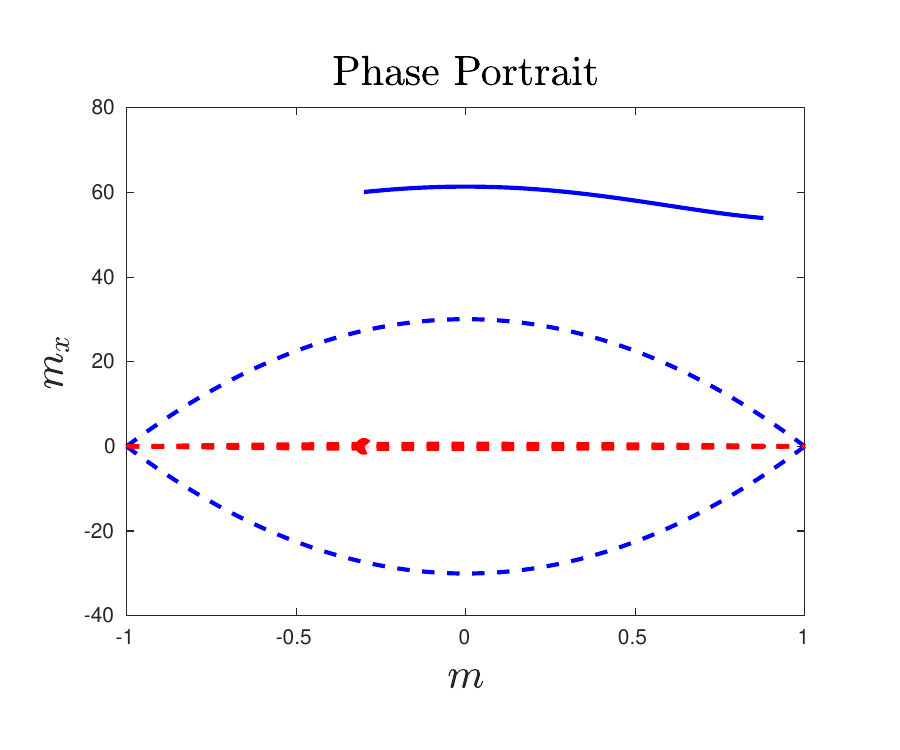}
        \caption{}
    \end{subfigure}
    \begin{subfigure}[c]{0.3\textwidth}
        \centering
         \includegraphics[width=\textwidth]{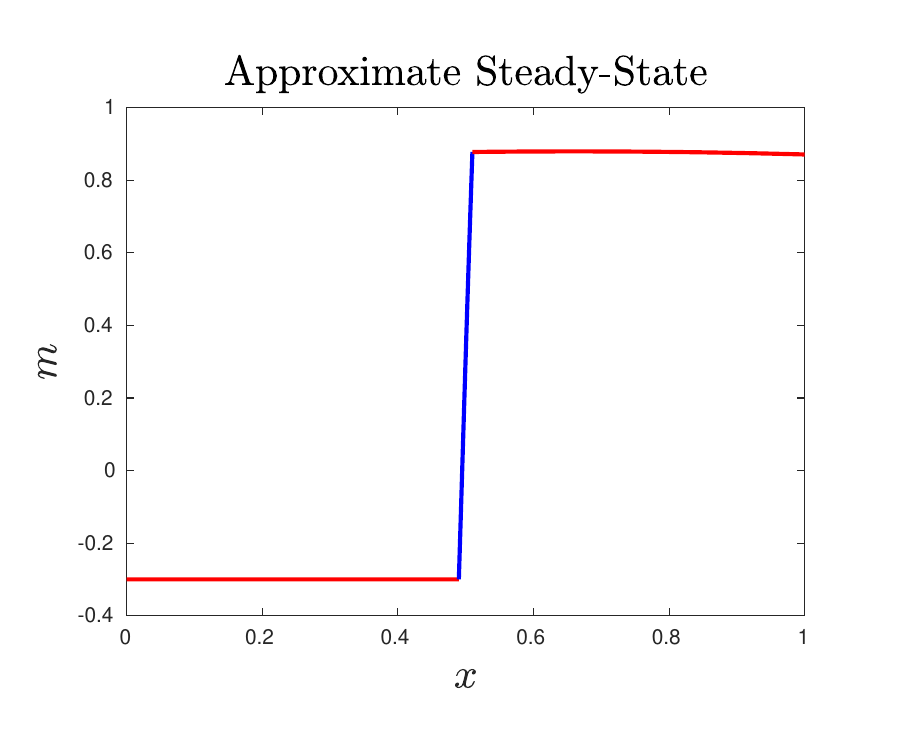}
        \caption{}
    \end{subfigure}
    \begin{subfigure}[c]{0.3\textwidth}
        \centering
         \includegraphics[width=\textwidth]{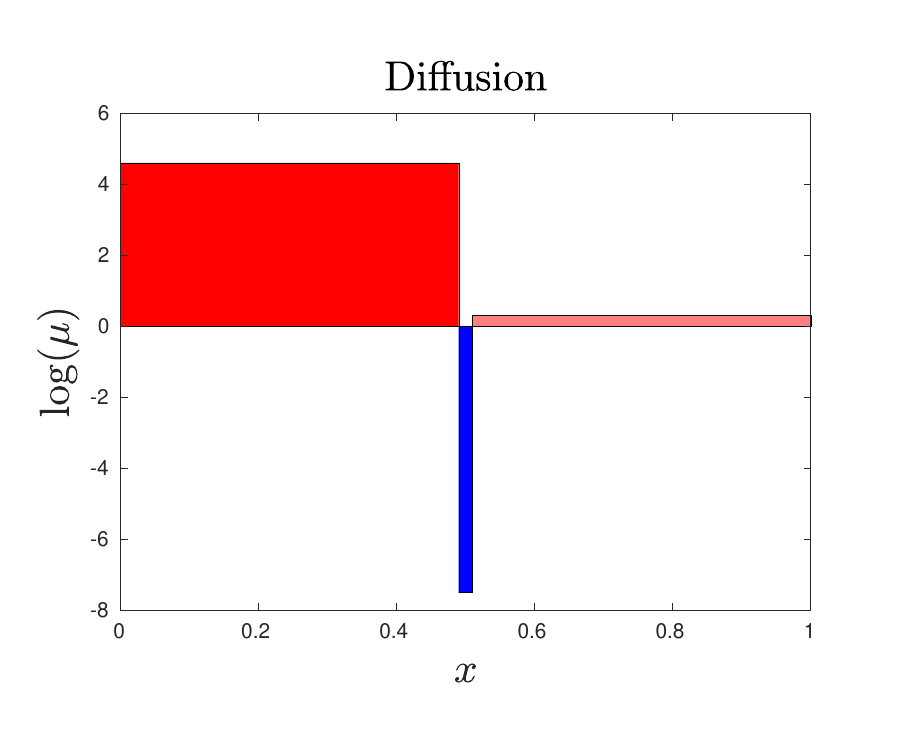}
        \caption{}
    \end{subfigure}
 \caption{(A) In dashed lines, the standing traveling wave delimiting the invariant region. The solid line denotes the trajectory employed. Red $\mu=100$, blue $\mu=5.537\times 10^{-4} $, light red $\mu=1.3647$. (B) Resulting approximate steady-state in the physical space. (C) Diffusion depicted in the physical space in a log scale.}\label{DIVexample}
\end{figure}

In Figure \ref{DIVexample}, we present an example of the procedure described and the diffusivity employed.

\subsection{Proof of Theorem \ref{TH2}}

Take $v\in L^\infty\left((0,1),[-1,1]\right)$ and let $S_v\in L^\infty\left((0,1),(-1,1)\right)$ be a piecewise constant function such that:
    $$\|v-S_v\|_{L^2}<\frac{\epsilon}{2}$$
    For $\alpha\in(-1,1)$ we can approximate a constant segment $w\equiv \alpha$ in $(0,L)$ for every $L>0$ by:
    \begin{equation}\label{pfeq1}
    \begin{cases}
          -\partial_{xx} m=\xi f(m)\\
     m=\alpha
    \end{cases}
    \end{equation}
    for $\xi$ small enough.
    
    Equation \eqref{pfeq1} has a solution in $L^\infty((0,1);[-1,1])$ since $m\equiv -1$ is a subsolution and $m\equiv 1$ is a supersolution. Therefore, we can consider a solution,  within such bounds, and so that the following estimate holds:
    $$|m_x|=\left|\int_0^L -\xi f(m)dx\right|\leq \xi L\|f\|_\infty$$
    and thus,
    $$\|m-\alpha\|_\infty=\left|\int_0^Lm_xdx\right|\leq \xi L^2\|f\|_\infty.$$

    The main problem is how to generate a steady-state  approximating a piecewise constant function made out of two components. We know that we can approximate a constant arc, and we could approximate several ones in different intervals. The remaining task is to find a steady state (defined in a small segment) whose Dirichlet and Neumann traces matches the two steady states generated with $\xi$ small. In this way, the steady state would make the transition from one segment to the other in a way that the $L^\infty$ bounds are respected.
    This is where the exact controllability of Lemma \ref{Col1} is of crucial importance.    
Using the controllability, we can approximate 
$$S_v:=\begin{cases}
        \alpha \quad 0\leq x\leq L,\\
        \beta\quad L<x\leq 1,
       \end{cases}
$$
with $\alpha,\beta\in(-1,1)$, by a steady-state considering the solutions of:
    \begin{equation}\label{pfeq2}
    \begin{cases}
          -\partial_{xx} m_1=\xi_1 f(m_1)&\quad x\in(0,L-\frac{\delta}{2}),\\
     m_1=\alpha &\quad x\in \{0,L-\frac{\delta}{2}\},
    \end{cases}
    \end{equation}
    and
        \begin{equation}\label{pfeq3}
    \begin{cases}
          -\partial_{xx} m_2=\xi_2f(m_2)&\quad x\in(L-\frac{\delta}{2},1),\\
     m_2=\beta &\quad x\in \{L-\frac{\delta}{2},1\}.
    \end{cases}
    \end{equation}
    We have to find a way to match the traces, in order to do so we formulate the following control problem     
     \begin{equation}\label{ODEcontr2}
        \begin{cases}
         \frac{d}{dx}
            \begin{pmatrix}
              m\\
              m_x
             \end{pmatrix}=\begin{pmatrix}
             m_x\\
             -\xi(x)f(m)
             \end{pmatrix},\qquad \xi(x)>0,\quad x\in(0,\delta),\\
             \vspace{0.1cm}\\
             \begin{pmatrix}
              m(0)\\
              m_x(0)
             \end{pmatrix}=\begin{pmatrix}
             m_1(L-\frac{\delta}{2})\\
             \frac{d}{dx}m_1(L-\frac{\delta}{2})
             \end{pmatrix},\qquad \begin{pmatrix}
              m(\delta)\\
              m_x(\delta)
             \end{pmatrix}=\begin{pmatrix}
             m_2(L+\frac{\delta}{2})\\
             \frac{d}{dx}m_2(L+\frac{\delta}{2})
             \end{pmatrix}.
             \end{cases}
\end{equation}
 The exact controllability of \eqref{ODEcontr2} with a trajectory inside $R$ is guaranteed by Lemma \ref{Col1} that we will prove subsequently (see subsection \ref{nondivlemma}). 
 Then 
\begin{align*}
 \|S_v-m\|_{L^2(L-\frac{\delta}{2},L+\frac{\delta}{2})}^2=\int_{L-\frac{\delta}{2}}^L |m-\alpha|^2dx&+\int_L^{L+\frac{\delta}{2}}|m-\beta|^2dx\leq 4\delta.
\end{align*}
Choosing $\delta^{1/2}=\epsilon/4$, and $\xi$ small enough in problems \eqref{pfeq2} and \eqref{pfeq3}, one has:
$$\| S_v-m\|_{L^2((0,1))}\leq \frac{\epsilon}{2}.$$

\subsubsection{Proof of Lemma \ref{Col1}}\label{nondivlemma}
\textcolor{white}{.}\newline
It suffices to understand the switching dynamics. 
Let us consider the initial data $(m_0,\partial_x m_0)$ and the target $(m_L,\partial_x m_L)$. We will start the trajectory in the initial data with a multiplicative control $\xi_1$ until the trajectory reaches its maximum velocity, at $m=0$, i.e at the point $(0,\partial_x \overline{m})$. At this point, we change the multiplicative control to a value $\xi_2$ so that the trajectory will eventually pass (not necessarily at $x=L$) through the target $(m_L,\partial_x m_L)$. We can ensure this by  the fact that both phases  employ a constant control, and hence the dynamics preserves the energy in each phase:
\begin{equation*}
 \begin{cases}
  \frac{1}{2}\partial_xm_0^2+\xi_1F(m_0)=\frac{1}{2}\partial_x \overline{m}^2+\xi_1F(0)\\
    \frac{1}{2}\partial_x \overline{m}^2+\xi_2F(0)=\frac{1}{2}\partial_xm_L^2+\xi_2F(m_L).
 \end{cases}
\end{equation*}
Note that $\xi_2$ as a function of $\xi_1$ is monotonous, a decrease (increase) on $\xi_1$ implies a decrease  (increase) in $\xi_2$.
Indeed, 
$$ \xi_2(\xi_1)=\frac{\partial_x m_0^2-\partial_xm_L^2}{2(F(m_L)-F(0))}+\frac{F(m_0)-F(0)}{F(m_L)-F(0)}\xi_1.$$
The constraints require  the multiplicative control $\xi$ to be positive. Hence, we have to impose the condition
$$\xi_1>-\frac{1}{2}\frac{\partial_x m_0^2-\partial_xm_L^2}{F(m_0)-F(0)}.$$
Up to now, we only know that there is a controlled trajectory reaching the target at a certain $x$, not necessarily at $x=L$. For ensuring this, we will make use of Claim \ref{Period}.

The period of an orbit starting at $(m_{\max},0))$ is a decreasing function of $\xi$. Hence, one may find, $\xi_1$ large enough so that the length required to achieve the target $(m_L,\partial_x m_L)$, $\mathcal{L}(\xi_1)$ is smaller than $L$, in fact
$$ \lim_{\xi_1\to +\infty} \mathcal{L}(\xi_1)=0.$$
Using that the period is a continuously decreasing function of $\xi$, one can find an orbit passing through the target with the period $X(\xi_3)$ being so that:
$$\frac{L-\mathcal{L}(\xi_1)}{X(\xi_3)}\in \mathbb{N}\backslash \{0\}.$$

\begin{figure}
\captionsetup[subfigure]{justification=centering}
\centering
    \begin{subfigure}[c]{0.45\textwidth}
        \centering
         \includegraphics[width=\textwidth]{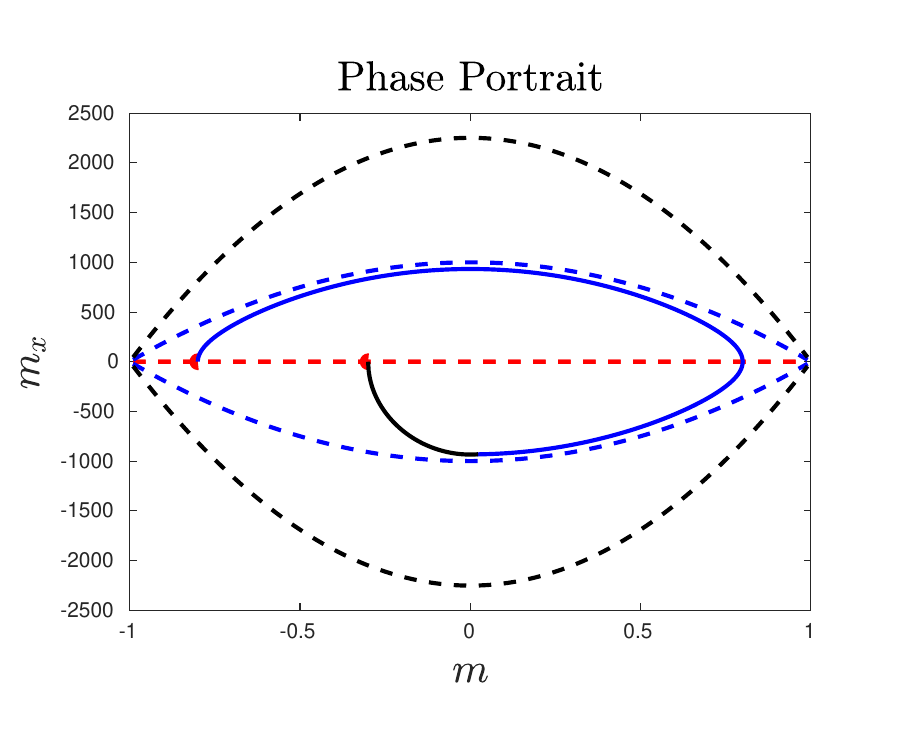}
        \caption{}
    \end{subfigure}
    \begin{subfigure}[c]{0.45\textwidth}
        \centering
         \includegraphics[width=\textwidth]{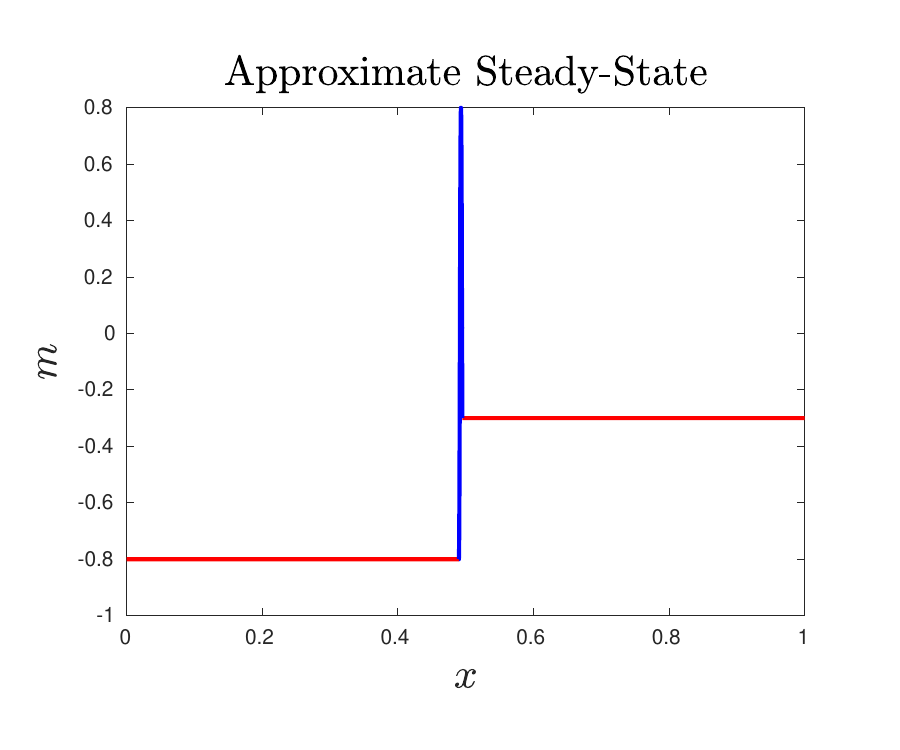}
        \caption{}
    \end{subfigure}

 \caption{(A) The dashed lines correspond to the traveling wave solutions limiting the invariant regions for different diffusivities, the filled curve is the controlled trajectory. Red $\mu=100$, blue $\mu=5\times 10^{-6}$, black $\mu=2\times 10^{-5}$. (B) representation of the approximate steady-state in the physical space.}\label{nDIVexample2}
\end{figure}

\subsection{Proof of Theorem \ref{TH1}}\label{proofTH1}
\begin{figure}
\captionsetup[subfigure]{justification=centering}
\centering
    \begin{subfigure}[c]{0.3\textwidth}
        \centering
         \includegraphics[width=\textwidth]{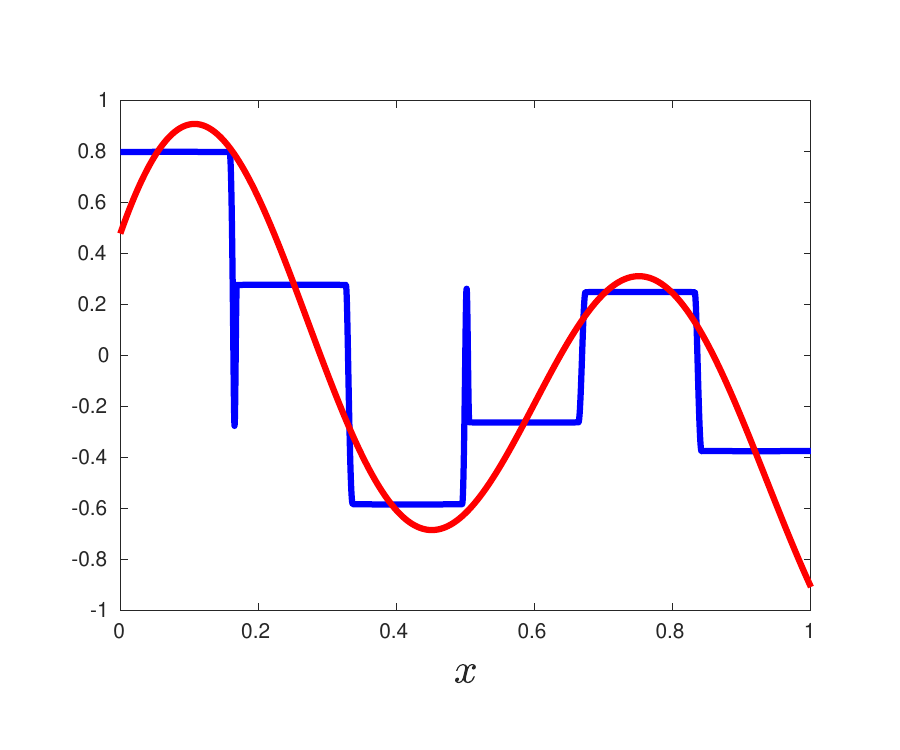}
        \caption{}
    \end{subfigure}
        \begin{subfigure}[c]{0.3\textwidth}
        \centering
         \includegraphics[width=\textwidth]{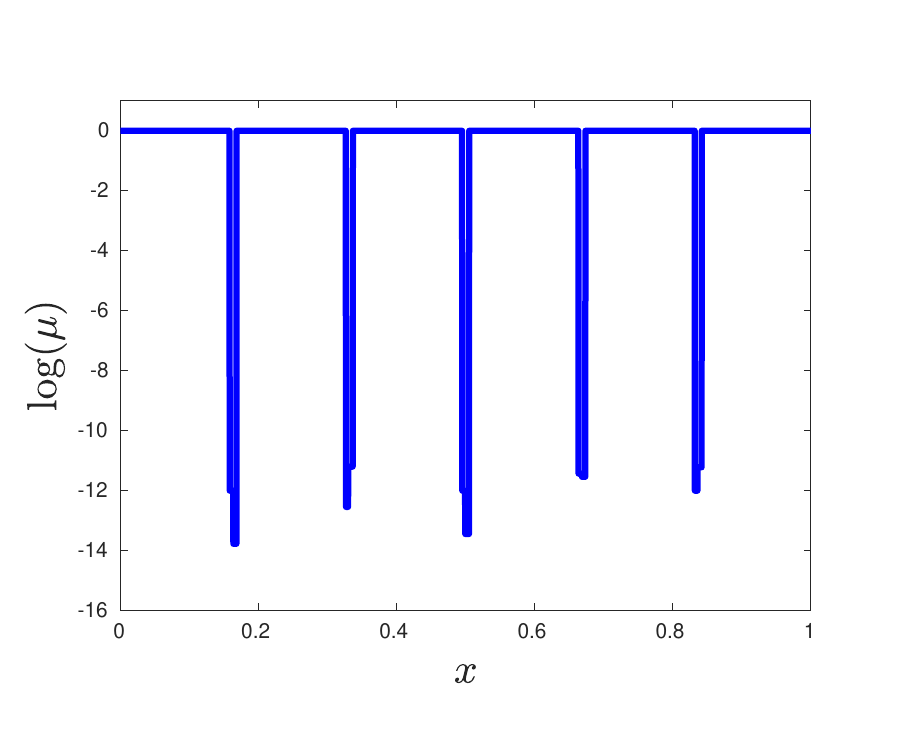}
        \caption{}
    \end{subfigure}
        \begin{subfigure}[c]{0.3\textwidth}
        \centering
         \includegraphics[width=\textwidth]{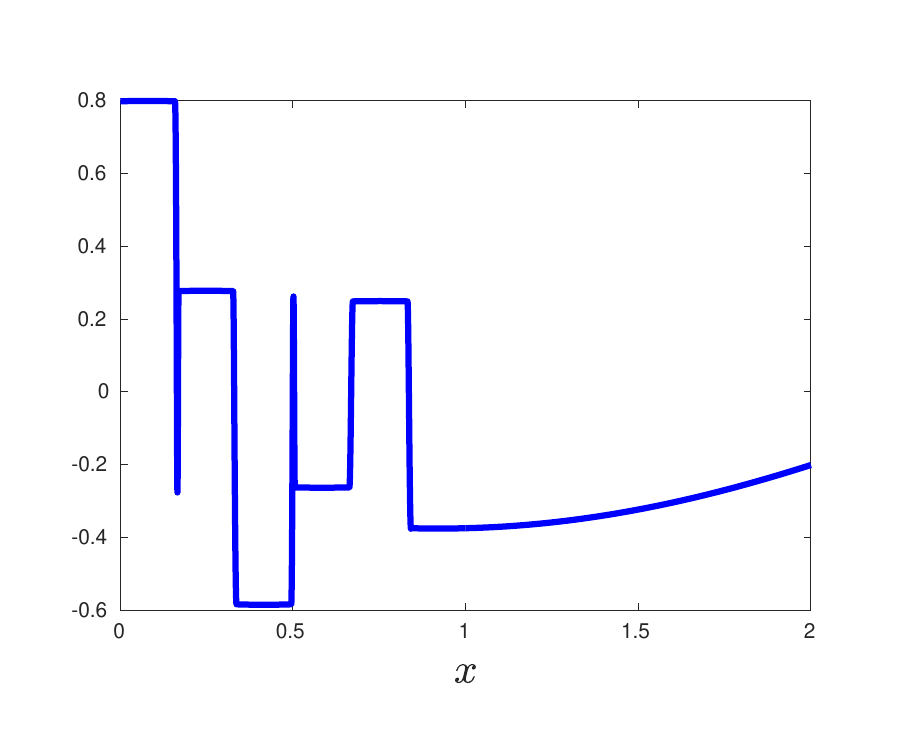}
        \caption{}\label{path}
    \end{subfigure}
   
 \caption{(A) In red, a random target function constructed randomizing the first three Fourier coefficients. In blue, an approximate steady-state. (B) The corresponding logarithm of the diffusion giving rise to the approximate steady-sate. (C) Extension of the approximate steady-state in $(0,2)$ with constant diffusion $\mu(x)$ in $(1,2)$.}\label{apprSS}
 \end{figure}

\textcolor{white}{.}\newline
Once the approximate target steady-state is built, see Figure \ref{apprSS}, we fix an arbitrary initial datum $m_0$ fulfilling the state-constraints, and we proceed as follows:
\begin{enumerate}
 \item  Set a constant very large diffusivity $\mu(x)\equiv\mu$  and $a=0$. Then, $m\equiv 0$ is stable, and the solution exponentially converges to it. Once we are close enough to $m=\equiv 0$, we apply local controllability from the boundary to control exactly at $m\equiv 0$. 
 \item We construct a path connecting $m\equiv 0$ to the approximate steady-state, and we employ the staircase method.
\end{enumerate}

The construction of the path is the following: 

Consider the approximate steady-state and its diffusion $(m_\epsilon,\mu_\epsilon)$. For the seek of simplicity we consider an approximate steady-state fulfilling the following condition on the boundary:
\begin{equation}\label{inv}
 \frac{1}{2}\partial_x m_\epsilon(1)<\frac{1}{\mu_\epsilon(1)}F(m_\epsilon(1)).
\end{equation}

Then, we can extend the steady-state to $(0,2)$ (see Figure \ref{path}) by employing the ODE dynamics extending the diffusion by a constant with continuity:
$$ \mu_E(x):=\begin{cases}
            \mu_\epsilon(x) &\quad x\in[0,1]\\
            \mu_\epsilon(1) &\quad x\in(1,2]
           \end{cases}$$

Since  \eqref{inv} holds the extended steady-state $m_E$ also satisfies the bounds:
$$ -1<m(x)<1 \quad x\in(0,2).$$

Now by shifting smoothly along $m_E$ we obtain a $s$-dependent family defined by 
\begin{equation*}
 \begin{cases}
  m_s(x)=m_E(x+s)\mathbb{1}_{(0,1)}\\
  \mu_s(x)=\mu_E(x+s)\mathbb{1}_{(0,1)}\\
  a_1=m_E(s)\quad a_2(s)=m_E(1+s)
 \end{cases}
\end{equation*}
which is continuous with respect to the $L^2$-norm,  connecting the approximate steady-state with a steady-state with homogeneous diffusion.

Then we can employ the techniques in \cite{PoucholTrelatZuazua,Domenec-Zuazua,drift} to find a path that connects with $m\equiv 0$.

Since the admissible path can be taken in the opposite way, we obtain the desired result.

The employed diffusivity coefficients are in $BV$. But thanks to the results of \cite{le2007carleman}, the staircase method can also be extended in such setting.

\section{Open problems and perspectives}\label{Open}

\textcolor{white}{.}\\
\textbf{1. Other nonlinearities.}
The fact of being able to orbit around $(0,0)$ in the ODE dynamics has been crucial for approximating any steady-state. It would be interesting to analyse whether similar results can be achieved for other bistable nonlinearities with $F(1)\neq 0$. 

Our proof relies on the existence of  an oscillator  and two trajectories arriving to and going away from $(1,0)$ and $(-1,0)$. In the case where $F(1)=0$ the invariant region fulfils both conditions. But for the case of $F(1)>0$, for instance, we would have to adapt the proof considering the invariant region containing periodic orbits with the trajectory associated to $E=F(1)$.

However, one can see that this type of result cannot be true for general nonlinearities. For instance, for monostable nonlinearities such as $f(s)=(s+1)(1-s)$ similar results cannot be achieved. Steady-states with a monostable nonlinearity reach their  minimum values always on the boundary and never in the interior (except for the trivial stationary solutions $m\equiv -1$ and $m\equiv 1$). Accordingly, a target with a minimum in the interior cannot be approximated. 

\textcolor{white}{.}\\
\textbf{2. Several space dimensions.}
In our $1-d$ constructions sharp transitions on the diffusivity have played a major role. Of course, the multi-dimensional setting presents a much higher geometric complexity. It is therefore unclear how the transitions on the diffusivity would need to be designed.
  
The connectivity result for  the set of steady-states can be extended to the multi-dimensional setting. \\

\noindent \textbf{3. Minimal controllability time.}
The analysis and numerical investigation of the minimal controllability time is a widely open subject in the present context.

\textcolor{white}{.}\\
\noindent \textbf{4. Other type of equations.}
Similar problems are of interest for   the Fokker-Planck diffusivity law \cite{van2005applicability}
\begin{equation}\label{Fock}
 \begin{cases}
  \partial_t m-\partial_{xx}(\mu(x,t)m)=f(m)\\
  m=a
 \end{cases}
\end{equation}

Nonlinear models such as the porous medium equation would also be interest:
\begin{equation}
 \begin{cases}
  \partial_t m-\partial_{x}(\mu(x,t)\partial_x(m^\alpha))=f(m)\\
  m=a
 \end{cases}
\end{equation}
Its analysis would require a more technical of the phase plane behaviour of trajectories.

Note that  for this model, due to the degeneration of the diffusivity when $m \sim 0$, not even the boundary controllability property without constraints is known (see \cite{coron2014global,geshkovski2019null}).
\\

\noindent \textbf{5. ODE Controllability problems.}
The analysis of this type of problems leads to interesting questions on ODE control as well. For instance, in the context of equation \eqref{Fock} 
 the following ODE control problem arises:
Setting $\mu(x,t)=\exp(g(x,t))$ to guarantee the positivity of the diffusion, and considering the ODE
\begin{equation}\label{ProposedControl}
 \frac{d}{dx}\begin{pmatrix}
  m\\
  m_x\\
  g\\
  g_x
 \end{pmatrix}
=\begin{pmatrix}
  m_x\\
  -e^{-g}f(m)-2m_xg_x-g_x^2m\\
  g_x\\
  0
 \end{pmatrix}
 +u(x)\begin{pmatrix}
  0\\
  m\\
  0\\
  1
 \end{pmatrix}
\end{equation}
one would need to understand its control properties with $u\in L^\infty((0,1))$ as control. 

For the PDE controllability result to hold,  it would be enough to prove the controllability for any $(m(0),m_x(0))\in R$ to any target function in $R$ while keeping the whole trajectory in $R$.

The possible limitations that could arise in the context of the control problem \eqref{ProposedControl} would  be inherited by the steady-states that the diffusion model can reach.

This task seems challenging  since \eqref{ProposedControl} is an affine control problem with drift, and  a global constrained control result is needed.
\\

\noindent \textbf{6. Efficient numerical methods.}
The development of efficient numerical solvers poses multiple difficulties, independent of the numerical scheme employed to discretise the PDE. Normally one should adopt an optimal control strategy. This would lead to non convex optimisation problems. Efficient methods would be required to deal, on one hand, with the complex pattern that the nearly optimal diffusion coefficients would adopt, and, on the other, with the necessity of fulfilling the constraints on the state and boundary controls.  
\textcolor{white}{.}\\

\noindent \textbf{7. Semilinear elliptic systems for pattern formation.}
One may consider a class of semilinear evolution equations such as
\begin{equation}\label{eqq}
 \begin{cases}
  \partial_t m^{(1)}-\mu_1\partial_{xx} m^{(1)}=a(x)f(m^{(1)})+b(x)m^{(2)}\quad & (x,t)\in(0,1)\times \mathbb{R}^+,\\
   \partial_t m^{(2)} -\mu_2\partial_{xx} m^{(2)}=c(x)m^{(2)}+d(x)g(m^{(2)})& (x,t)\in (0,1)\times \mathbb{R}^+,\\
   m^{(1)}=0,m^{(2)}=0 & (x,t)\in \{0,1\}\times \mathbb{R}^+,\\
      m^{(1)}=m^{(1)}_0,m^{(2)}=m^{(2)}_0 & (x,t)\in (0,1)\times \{0\},
 \end{cases}
\end{equation}
with $f,g\in C^1$, $f(0)=g(0)=0$ and $f'(0)=g'(0)=1$.
Suppose that, for a specific open set of initial data $\mathcal{X}\subset L^2((0,1))\times L^2((0,1))$, the omega limit of $\omega(\mathcal{X})$ for the dynamics \eqref{eqq} is a stable solution of the elliptic system
\begin{equation}\label{eqqq}
 \begin{cases}
  -\mu_1\partial_{xx} m^{(1)}=a(x)f(m^{(1)})+b(x)m^{(2)}\quad &  x\in(0,1),\\
    -\mu_2\partial_{xx} m^{(2)}=c(x)m^{(2)}+d(x)g(m^{(2)})\quad &  x\in(0,1),\\
       m^{(1)}=0,m^{(2)}=0 & x\in \{0,1\},
 \end{cases}
\end{equation}
This, of course, would be a subject of investigation. But, even if that were true, the question of determining the class of configurations that could be approximated by the steady-states of
\eqref{eqqq} would be also worth considering.

\section*{Acknowledgments}
The authors thank Dario Pighin for his valuable comments and discussions.

D. Ruiz-Balet was funded by the UK Engineering and Physical Sciences Research Council (EPSRC) grant EP/T024429/1. The second author has been funded by the Alexander von Humboldt-Professorship program, the Transregio 154 Project ``Mathematical Modelling, Simulation and Optimization Using the Example of Gas Networks" of the DFG, the ModConFlex Marie Curie Action, HORIZON-MSCA-2021-$d$N-01, AFOSR Proposal 24IOE027, grants PID2020-112617GB-C22 and TED2021-131390B-I00 of MINECO (Spain), and by the Madrid Government -- UAM Agreement for the Excellence of the University Research Staff in the context of the V PRICIT (Regional Programme of Research and Technological Innovation). 

\bibliography{CBIBfull29m.bib}{}
\bibliographystyle{amsplain}
\end{document}